\documentclass[12pt]{amsart}

  \usepackage{graphicx}
  \usepackage{color}
  \usepackage{boxedminipage}
  \usepackage{pictexwd,dcpic}
  \usepackage{latexsym,amscd,amssymb,epsfig,a4wide}

  \newcommand{\al }{\alpha }
  \newcommand{\cC }{\mathcal{C}}

  \newcommand{\Cm }{C}
  \newcommand{\cm }{c}
  
  \newcommand{\wo }{\le _R}
  \newcommand{\swo }{<_R}
  \newcommand{\Wgfun }{E}

  \newcommand{\Homsto }[2][\Wg (\cC )]{\mathrm{Hom}(#1,#2)}
  \newcommand{\id }{\mathrm{id}}
  \newcommand{\lspan }{\mathrm{span}}
  \newcommand{\ndN }{\mathbb{N}}
  
  \newcommand{\ndR }{\mathbb{R}}
  \newcommand{\ndZ }{\mathbb{Z}}
  \newcommand{\arrg }{\mathcal{A}}
  \newcommand{\Perm }[1]{\mathfrak{S}_{#1}}
  \newcommand{\rer }[1]{(R\re)^{#1}}
  \newcommand{\rfl }{\rho }
  
  \newcommand{\Ob }{\mathrm{Ob}}
  
  \newcommand{\re }{^\mathrm{re}}
  \newcommand{\rsC }{\mathcal{R}}
  \newcommand{\s }{\sigma }

  \newcommand{\Wg }{\mathcal{W}}
  \newcommand{\fac}{F}
  \DeclareMathOperator{\Aut}{Aut}
  
  \DeclareMathOperator{\Hom}{Hom}
  
  \def\RR{{\mathbb R}}

  \theoremstyle{plain}
    \newtheorem{theorem}{Theorem}[section]
    \newtheorem{proposition}[theorem]{Proposition}
    \newtheorem{lemma}[theorem]{Lemma}
    \newtheorem{corollary}[theorem]{Corollary}
    
  \theoremstyle{definition}
    \newtheorem{definition}[theorem]{Definition}
    \newtheorem{example}[theorem]{Example}

    \newtheorem{remark}[theorem]{Remark}

  \numberwithin{equation}{section}

  \newenvironment{note}[1][Note]
   {\bigskip\begin{center}\begin{boxedminipage}{4.5in}\setlength{\parindent}{1em}\noindent\textbf{#1.
}}
   {\end{boxedminipage}\end{center}\bigskip}

\begin{document}

  \title{Geometric Combinatorics of Weyl Groupoids}

  \author{Istv\'an Heckenberger}
  \thanks{I.H. was supported by the German
  Research Foundation via a Heisenberg fellowship}
  \address{Fachbereich Mathematik und Informatik\\
           Philipps-Universit\"at Marburg\\
           35032 Marburg, Germany}
  \email{heckenberger@mathematik.uni-marburg.de}

  \author{Volkmar~Welker}
  \email{welker@mathematik.uni-marburg.de}

  \begin{abstract}
    We extend properties of the weak order on finite
    Coxeter
    groups to Weyl groupoids admitting a finite root system.
    In particular, we determine the topological structure of intervals with
    respect to weak order, and show that the set of morphisms with fixed
    target object forms an ortho-complemented meet semilattice.
    We define the Coxeter complex of a Weyl groupoid with finite root system
    and show that it coincides with the triangulation of a sphere
    cut out by a simplicial hyperplane arrangement.
    As a consequence, one
    obtains an algebraic interpretation of many hyperplane arrangements
    that are not reflection arrangements.
  \end{abstract}

  \maketitle

  \section{Introduction}
  \label{sec:intro}

  Finite crystallographic Coxeter groups, also known as finite Weyl groups,
  play a prominent role in many branches of mathematics like combinatorics,
  Lie theory, number theory, and geometry. In the late sixties
  V.~Kac and R.\,V.~Moody (see \cite{b-Kac90}) discovered independently
  a class of infinite dimensional Lie algebras. In their approach the Weyl
  group is defined in terms of a generalized Cartan matrix. Later in the
  seventies V.~Kac also introduced Lie superalgebras using even more general
  Cartan matrices \cite{a-Kac77},
  and observed that different Cartan matrices may give rise to isomorphic Lie
  superalgebras. S.~Khoroshkin and V.~Tolstoy \cite[p.\,77]{inp-KhorTol95}
  observed that
  the Weyl group symmetry of simple Lie algebras can be generalized to a
  Weyl groupoid symmetry of contragredient Lie superalgebras, without working 
  out the details. Independently, Weyl groupoids turned out to be the main tool
  for the study of finiteness properties of Nichols algebras
  \cite{a-AndrSchn02} over groups. 

  Motivated by these developments, an axiomatic study
  of Weyl groupoids was initiated
  by H.~Yamane and the first author \cite{a-HeckYam08}.
  The theory was further extended by a series of papers of M.~Cuntz and the
  first author, and a satisfactory classification result of finite Weyl
  groupoids of rank two and three was achieved \cite{a-CH09b,p-CH09b}.
  Interestingly, not all finite Weyl groupoids obtained via the classification
  are related to known Nichols algebras. A possible explanation could be
  the existence
  of an additional axiom which holds for the Weyl groupoid of any Nichols
  algebra. However, no such axiom was found yet, and a more systematic
  study is needed to find some clue.

  Coxeter groups, in particular Weyl groups, are a source of important classes
  of examples for simplicial hyperplane arrangements (see for example the
  seminal work of P.~Deligne \cite{a-Deligne72}).
  Roughly speaking, a
  simplicial hyperplane arrangement is a family of hyperplanes in a Euclidean
  space that cuts space into simplicial cones. However, most simplicial
  arrangements have no interpretation in terms of Coxeter groups. Therefore
  there is no canonical algebraic structure which hints toward a description 
  of the fundamental group of the complement of the complexification as
  described in \cite{a-Deligne72}. Also in general simplicial arrangements
  lack a relation to Lie algebras. It was observed in \cite{p-CH09b}
  that Weyl groupoids of rank three are related to simplicial arrangements in
  a real projective plane. Interestingly, the classification of
  such arrangements is not yet completed \cite{a-Gruenb09}.
  It was noted in \cite{p-CH09b} that most known exceptional
  arrangements, in particular the largest one,
  can be explained via Weyl groupoids.

  In this paper we analyze the structure of the Weyl groupoid related to
  parabolic subgroups and the weak order. Most of our results are known for
  Coxeter groups from the work of A.~Bj\"orner 
  (see \cite{a-Bj84}, \cite{inp-Bj84}, \cite{b-BjoeBre05}). 
  Our goal is to find an appropriate
  generalization. For the proofs either a careful adaption of the classical
  proofs is required or the lack of group structure forces new proofs which
  in some cases seem to be simpler than the usual ones.

  The weak order is defined using the length function on the Weyl
  groupoid. It proved its relevance for Coxeter groups,
  and it also has an interpretation for Nichols algebras in terms of right
  coideal subalgebras \cite{p-HeckSchn09}. We work out an example
  (Example \ref{ex:bruhat}) which
  shows that the weak order on a Weyl groupoid may have significantly different
  properties than the one on a Coxeter group. As a consequence, our results
  cover a much wider class of partially ordered sets and simplicial
  arrangements than the classical ones. We investigate
  longest elements of parabolic subgroupoids, and show in
  Proposition~\ref{pr:extdesc} that the poset they define is isomorphic to the
  poset of subsets of the set of simple reflections. In
  Theorem~\ref{th:semilattice} we prove that the set of morphisms with
  fixed target object is a meet semilattice. It is worthwhile to mention that
  this result is usually proved using the exchange condition, which is not
  available for Weyl groupoids \cite{a-HeckYam08}.
  For our proof we take advantage of our knowledge on longest elements.
  In addition with Theorem~\ref{th:Ju} we find a formula involving the
  letters of the meet of two words in the weak order.
  With Theorem~\ref{th:interval} we clarify the topological structure of
  intervals in weak order, and in Theorem~\ref{th:orthocomplemented}
  it is shown that the set of morphisms with fixed target object is
  ortho-complemented.

  In Section~\ref{sec:coxeter} we give two different definitions of the
  Coxeter complex associated to a fixed object of a Weyl groupoid. From
  one of the definitions it is immediate that the Coxeter complex is
  simplicial, and the other one shows that it comes from a hyperplane
  arrangement. We prove in Corollary ~\ref{co:isomorphism} that the two
  definitions yield isomorphic complexes, and hence the Coxeter complex is a
  simplicial complex which can be seen as the complex induced by a
  simplicial hyperplane
  arrangement on the unit sphere.

  \section{Basic Concepts}
  \label{sec:basics}

  \subsection{Weyl groupoids}

  We mainly follow the notation in \cite{a-CH09a,a-CH09b}.
  The fundaments of the general theory have been developed in \cite{a-HeckYam08}.
  Let us start by recalling the main definitions.

  Let $I$ and $A$ be finite sets with $A\not=\emptyset $. Let
  $\{\al _i\,|\,i\in I\}$ be the standard basis of $\ndZ ^I$.
  For all $i\in I$
  let $\rfl _i : A \to A$ be a map,
  and for all $a\in A$ let $\Cm ^a=(\cm ^a_{jk})_{j,k \in I}$
  be a generalized Cartan matrix
  in the sense of \cite[\S 1.1]{b-Kac90}, where $\cm ^a_{j k}\in \ndZ $ for
  all $j,k\in I$.
  The quadruple
  \[ \cC = \cC (I,A,(\rfl _i)_{i \in I}, (\Cm ^a)_{a \in A})\]
  is called a \textit{Cartan scheme} if
  \begin{enumerate}
    \item[(C1)] $\rfl _i^2 = \id$ for all $i \in I$,
    \item[(C2)] $\cm ^a_{ij} = \cm ^{\rfl _i(a)}_{ij}$ for all $a\in A$ and
          $i,j\in I$.
  \end{enumerate}

  Let $\cC = \cC (I,A,(\rfl _i)_{i \in I}, (\Cm ^a)_{a \in A})$ be a
  Cartan scheme. For all $i \in I$ and $a \in A$ define $\s _i^a \in
  \Aut(\ndZ ^I)$ by
  \begin{align}
    \s _i^a (\al _j) = \al _j - \cm _{ij}^a \al _i \qquad
    \text{for all $j \in I$.}
    \label{eq:sia}
  \end{align}
  Then $\s _i^a$ is a reflection in the sense of
  \cite[Ch.\,V,\,\S\,2]{b-BourLie4-6}.
  The \textit{Weyl groupoid of} $\cC $
  is the category $\Wg (\cC )$ such that $\Ob (\Wg (\cC ))=A$ and
  the morphisms are compositions of maps
  $\s _i^a$ with $i\in I$ and $a\in A$,
  where $\s _i^a$ is considered as an element in $\Hom (a,\rfl _i(a))$.
  The category $\Wg (\cC )$ is a groupoid.
  The set of morphisms of $\Wg (\cC )$ is also denoted by $\Wg (\cC )$,
  and we use the notation
  \[ \Homsto{a}=\mathop{\cup }_{b\in A}\Hom (b,a) \quad
  \text{(disjoint union)}. \]

  \begin{example} \label{ex:classical}
    Let $(W,S)$ be a Coxeter system for a crystallographic Coxeter group $W$. 
    Then $(W,S)$ can be seen as a Weyl groupoid $\Wg(\cC)$ with a single object 
    $a$ and $\Hom(a,a) = \langle S \rangle = W$ with Cartan scheme $\cC = 
    \cC(\{1, \ldots, |S|\},\{ a\}, (\rfl _i=\id )_{i=1,...,|S|}, (C^a))$
    where $C^a$ is the usual Cartan matrix of $W$.
  \end{example}

  For notational convenience we will often neglect upper indices referring to
  elements of $A$ if they are uniquely determined by the context. For example,
  the morphism
  \[ \s _{i_1}^{\rfl _{i_2}\cdots \rfl _{i_k}(a)}
  \cdots \s_{i_{k-1}}^{\rfl _{i_k}(a)}\s _{i_k}^a\in \Hom (a,b),
  \quad \text{where $k\in \ndN _0$, $i_1,\dots,i_k\in I$, and
  $b=\rfl _{i_1}\cdots \rfl _{i_k}(a)$,} \]
  will be denoted by $\s _{i_1}\cdots \s _{i_k}^a$ or by
  $\id^b\s _{i_1}\cdots \s_{i_k}$.
  The cardinality of $I$ is termed the \textit{rank of} $\Wg (\cC )$.
  A Cartan scheme is called \textit{connected} if its Weyl groupoid
  is connected, that is, if for all $a,b\in A$ there exists $w\in \Hom (a,b)$.
  The Cartan scheme is called \textit{simply connected},
  if for all $a,b \in A$ the set $\Hom (a,b)$ consists of at most one element.

  Let $\cC $ be a Cartan scheme. For all $a\in A$ let
  \[ \rer a=\{ \id^a \s _{i_1}\cdots \s_{i_k}(\al _j)\,|\,
  k\in \ndN _0,\,i_1,\dots,i_k,j\in I\}\subseteq \ndZ ^I.\]
  The elements of the set $\rer a$ are called \textit{real roots} (at $a$) --
  this notion is adopted from \cite[\S 5.1]{b-Kac90}.
  The pair $(\cC ,(\rer a)_{a\in A})$ is denoted by $\rsC \re (\cC )$.
  A real root $\al \in \rer a$, where $a\in A$, is called \textit{positive}
  (resp.\ \textit{negative})
  if $\al \in \ndN _0^I$ (resp.\ $\al \in -\ndN _0^I$).
  In contrast to real roots associated to a single generalized Cartan matrix
  (e.g. Example \ref{ex:classical}),
  $\rer a$ may contain elements which are neither positive nor negative. A good
  general theory can be obtained if $\rer a$ satisfies additional properties.

  Let $\cC =\cC (I,A,(\rfl _i)_{i\in I},(\Cm ^a)_{a\in A})$ be a Cartan
  scheme. For all $a\in A$ let $R^a\subseteq \ndZ ^I$, and define
  $m_{i,j}^a= |R^a \cap (\ndN_0 \al _i + \ndN_0 \al _j)|$ for all $i,j\in
  I$ and $a\in A$. One says that
  \[ \rsC = \rsC (\cC , (R^a)_{a\in A}) \]
  is a \textit{root system of type} $\cC $, if it satisfies the following
  axioms.
  \begin{enumerate}
    \item[(R1)]
      $R^a=R^a_+\cup - R^a_+$, where $R^a_+=R^a\cap \ndN_0^I$, for all
      $a\in A$.
    \item[(R2)]
      $R^a\cap \ndZ\al _i=\{\al _i,-\al _i\}$ for all $i\in I$, $a\in A$.
    \item[(R3)]
      $\s _i^a(R^a) = R^{\rfl _i(a)}$ for all $i\in I$, $a\in A$.
    \item[(R4)]
      If $i,j\in I$ and $a\in A$ such that $i\not=j$ and $m_{i,j}^a$ is
      finite, then
      $(\rfl _i\rfl _j)^{m_{i,j}^a}(a)=a$.
  \end{enumerate}

  \begin{example} \label{classicalroot}
    Let $(W,S)$ be a Coxeter system for a finite
    crystallographic Coxeter group $W$
    acting on some real vector space $V$ 
    seen as a Weyl groupoid as in Example \ref{ex:classical}. Then by
    \cite[p. 6]{b-Humphreys90} a root system of $W$ is a set of vectors $R$
    from $V$ such that :
    \begin{itemize}
      \item[(R1')] $R \cap \RR \al = \{ \al,-\al\}$ for all $\al \in R$.
      \item[(R2')] $\s R = R$ for all reflections $\s$ from $W$.
    \end{itemize}
    Clearly, (R1') implies (R2) and from the finiteness and
    the crystallographic condition we infer
    that (R2) implies (R1'). It is obvious that (R2') implies 
    (R3). Since any reflection is a product of simple reflections
    it follows that
    (R3) implies (R2'). Since our groupoid has only one object,
    Axiom (R4) is vacuous. 
    As a consequence \cite[p. 8]{b-Humphreys90} of (R1') and
    (R2') every set of positive roots contains a unique simple system. Then the
    definition of a simple system and the crystallographic condition
    imply (R1).
    Thus we have shown that for finite
    crystallographic Coxeter groups conditions
    (R1')-(R2') and (R1)-(R3) are equivalent.
  \end{example}

  Axioms (R2) and (R3) are always fulfilled for $\rsC \re $.
  A root system $\rsC $ is called \textit{finite} if for all $a\in A$ the
  set $R^a$ is finite. By \cite[Prop.\,2.12]{a-CH09a},
  if $\rsC $ is a finite root system
  of type $\cC $, then $\rsC =\rsC \re $, and hence $\rsC \re $ is a root
  system of type $\cC $ in that case.

  In \cite[Def.\,4.3]{a-CH09a} the concept of an \textit{irreducible}
  root system of type
  $\cC $ was defined. By \cite[Prop.\,4.6]{a-CH09a}, if $\cC $ is a connected
  Cartan scheme and $\rsC $ is a finite root system of type $\cC $,
  then $\rsC $ is irreducible if and only if for all $a\in A$ (or,
  equivalently, for some $a\in A$)
  the generalized Cartan matrix $C^a$ is indecomposable.

  Let $\cC =\cC (I,A,(\rfl _i)_{i\in I},(\Cm ^a)_{a\in A})$ be a Cartan
  scheme.
  Let $\Gamma $ be an undirected graph,
  such that the vertices of $\Gamma $ correspond to the elements of $A$.
  Assume that for all $i\in I$ and $a\in A$ with $\rfl _i(a)\not=a$
  there is precisely one edge between the vertices $a$ and $\rfl _i(a)$
  with label $i$,
  and all edges of $\Gamma $ are given in this way.
  The graph $\Gamma $ is called the \textit{object change diagram} of $\cC $.

  Now we introduce parabolic subgroupoids which will play a crucial role
  in the sequel.

  \begin{definition}
    Let $\cC =\cC (I,A,(\rfl _i)_{i\in I},(\Cm ^a)_{a\in A})$ be a Cartan
    scheme and let $J\subseteq I$.
    The {\em parabolic subgroupoid} $\Wg _J(\cC )$
    is the smallest subgroupoid of $\Wg (\cC )$ which contains all objects
    of $\Wg (\cC )$ and all morphisms $\s _j^a$ with $j\in J$ and $a\in A$.
  \end{definition}

  In general, parabolic subgroupoids are not connected, even if $\cC $
  is connected.

  The most important tools
  for the study of the weak order in the next section
  will be the length
  functions of the parabolic subgroupoids $\Wg _J(\cC )$ of $\Wg (\cC )$,
  where $J\subseteq I$.
  For all $J\subseteq I$ let $\ell _J:\Wg _J(\cC )\to \ndN _0$ such
  that
  \begin{align}
    \ell _J(w)=\min \{k\in \ndN _0\,|\,w=\s _{i_1}\cdots \s _{i_k}^a,
    i_1,\dots,i_k\in J\}
    \label{eq:ellJ}
  \end{align}
  for all $a,b\in A$ and $w\in \Hom (a,b)$. For $J=I$ this is the adaption of
  the usual length
  function from classical Coxeter groups to Weyl groupoids defined in \cite{a-HeckYam08}.
  We write $\ell (w)$ instead of $\ell _I(w)$. For $w \in \Wg(\cC)$
  we say that $w = \s_{i_1}
  \cdots \s_{i_k}$ is a {\it reduced decomposition} of $w$ if $k = \ell(w)$.  

  The length
  function on Weyl groupoids has similar properties as the usual length
  function on Coxeter groups, see \cite{a-HeckYam08}.
  In particular the following holds.

  \begin{lemma}[Lemma 8(iii) \cite{a-HeckYam08}] \label{le:length}
    Let $a,b\in A$ and $w\in \Hom (a,b)$. Then
    \[ \ell (w)=|\{\al \in R^a_+\,|\,w(\al )\in -R^b_+\}|. \]
  \end{lemma}

  \begin{lemma}[Corollary 3 \cite{a-HeckYam08}] \label{le:posi}
    Let $a,b\in A$, $w\in \Hom (a,b)$, and $i\in I$. Then $\ell (w\s _i)=\ell
    (w)-1$ if and only if $w(\al _i)\in -R^b_+$. Equivalently, $\ell (w\s
    _i)=\ell (w)+1$ if and only if $w(\al _i)\in R^b_+$.
  \end{lemma}

  Before we proceed with studying the length function itself we clarify the
  structure of the set of subsets $J \subseteq I$ for which 
  $w \in \Hom(a,b)$ is also a morphism in $\Wg_J (\cC)$.

  \begin{proposition} \label{pr:content}  
    Let $w \in \Hom(a,b)$. If $w = \s_{i_1} \cdots \s_{i_k}^a$ is a
    reduced decomposition of $w$ and
    $w = \s_{j_1} \cdots \s_{j_l}^a$ is another decomposition, where $k,l\in
    \ndN _0$ and $i_1,\dots,i_k,j_1,\dots,j_l\in I$, 
    then as sets $$\{ i_1, \ldots, i_k \} \subseteq \{ j_1, \ldots, j_l \}.$$
    In particular, if $k=l$ then
    $\{ i_1, \ldots, i_k \}=\{ j_1, \ldots, j_k \}$.
  \end{proposition}

  \begin{proof}
    Set $J := \{ i_1, \ldots, i_k \}$ and $J' = \{ j_1, \ldots, j_l \}$. 
    Assume that $J \not\subseteq J'$.
    Let $m\in \{1,\dots,k\}$ such that $i_m\notin J'$ and $i_{m'}\in J'$ for
    all $m'<m$. Let $\al =\id^a\s _{i_k} \s_{i_{k-1}} \cdots
    \s_{i_{m+1}}(\alpha_{i_m})$.
    Then $\al \in R^a_+$ by the fact that $w = \s_{i_1} \cdots \s_{i_k}^a$ is
    a reduced decomposition and by Lemma~\ref{le:posi}. Moreover, 
    \begin{align}
      w(\al )=\s _{i_1}\cdots \s _{i_{m-1}}\s _{i_m}(\al _{i_m})
      =-\s _{i_1}\cdots \s _{i_{m-1}}(\al _{i_m})\in -\al _{i_m}+
      \lspan _\ndZ \{\al _j\,|\,j\in J'\}.
      \label{eq:wal}
    \end{align}
    Let $\al =\al '+\al ''$ with
    $\al '\in \lspan _{\ndN _0}\{\al _j\,|\,j\notin J'\}$ and
    $\al ''\in \lspan _{\ndN _0} \{\al _j\,|\,j\in J'\}$.
    Since $w\in \Wg _{J'}(\cC )$, we conclude that $w(\al )\in \al '
    +\lspan _\ndZ \{\al _j\,|\,j\in J'\}$. This is a contradiction to
    \eqref{eq:wal} since $i_m\notin J'$. Hence $J\subseteq J'$. 
  \end{proof}

  For all $a,b \in A$, $w \in \Hom(a,b)$ and reduced decompositions 
  $w = \s_{i_1} \cdots \s_{i_k}^a$ we set $J(w) := \{ i_1, \ldots, i_k \}$.
  By Proposition \ref{pr:content} this definition is independent of the
  chosen reduced decomposition. Moreover, for any subset $J\subseteq I$ and
  any $w\in \Wg _J(\cC )$ the reduced decompositions of $w$ are also contained
  in $\Wg _J(\cC )$.
  Observe also that $J(w)=J(w^{-1})$ for all
  $w\in \Wg (\cC )$ and that $J(uv)=J(u)\cup J(v)$
  for all $u,v\in \Wg (\cC )$ with $\ell (uv)=\ell (u)+\ell (v)$.

  \begin{corollary} \label{co:lJ=l}
    Let $J\subseteq I$. Then $\ell _J(w)=\ell (w)$ for all
    $w\in \Wg _J(\cC )$.
  \end{corollary} 

  \begin{proof} 
    If there is a decomposition of $w$ having only factors $\s_i$ with
    $i\in J$
    then by Proposition \ref{pr:content} all reduced decompositions have this
    property.
    The assertion follows.
  \end{proof}

  One can characterize $J(w)$ for any $w\in \Wg (\cC )$ in terms of roots.

  \begin{lemma}
    Let $a,b\in A$, $J\subseteq I$, and let $w\in \Hom (b,a)$. Then
    $J(w)\subseteq J$ if and only if
    $w(R^b_+)\subseteq R^a_+\cup \sum _{j\in J}\ndZ \al _j$.
    \label{le:J(w)char}
  \end{lemma}

  \begin{proof}
    The implication $\Rightarrow$ follows from the definition of simple
    reflections and from Axioms (R1), (R3).
    Assume now that
    $w(R^b_+)\subseteq R^a_+\cup \sum _{j\in J}\ndZ \al _j$ and that $J(w)\not
    \subseteq J$. Then $J(\s _i w)\not\subseteq J$ and
    $\s _i w(R^b_+)\subseteq R^{\rfl _i(a)}_+\cup
    \sum _{j\in J}\ndZ \al _j$ for all $i\in J$,
    and
    hence by multiplying $w$ from the left by an appropriate element of $\Wg
    _J(\cC )$ we may assume that $\ell (\s _jw)=\ell (w)+1$ for all
    $j\in J$. It follows that $w^{-1}(\al _j)\in R^b_+$ for all $j\in J$ by
    Lemma~\ref{le:posi}. Hence $w(R^b_+)\subseteq R^a_+$, and therefore
    $w=\id ^a$ by Lemma~\ref{le:length}. This is a contradiction
    to $J(w)\not\subseteq J$.
  \end{proof}

  Let $J\subseteq I$ and for all $a\in A$ let
  $\Cm '^a=(\cm '^a_{j k})_{j,k\in J}$. Then
  $\cC '=\cC '(J,A,(\rfl _j)_{j\in J},(\Cm '^a)_{a\in A})$ is a Cartan scheme.
  It is denoted by $\cC |_J$ and is called the \textit{restriction of} $\cC $
  to $J$. As noted in \cite[Sect.\,4]{a-CH09a}, if $\rsC \re (\cC )$ is a
  root system of type $\cC $, then $\rsC \re (\cC |_J)$ is a root system of
  type $\cC |_J$, and finiteness of $\rsC \re (\cC )$ implies finiteness of
  $\rsC \re (\cC |_J)$. We compare restrictions with parabolic subgroupoids.

  \begin{lemma}
    Let $J\subseteq I$, $a\in A$, $k\in \ndN _0$,
    and $i_1,\dots,i_k\in J$ such that $\s _{i_1}\cdots \s_{i_k}^a|_{\ndZ
    ^J}=\id^a|_{\ndZ ^J}$. Then $\s _{i_1}\cdots \s _{i_k}^a=\id^a$.
    \label{le:wellness}
  \end{lemma}

  \begin{proof}
    By assumption $\s _{i_1}\cdots \s _{i_k}^a(\al _j)=\al _j$ for all $j\in
    J$. Since $i_1,\ldots, i_k\in J$,  the definition of $\s _j^b$ for $j\in
    J$, $b\in A$ implies that $\s _{i_1}\cdots \s _{i_k}^a(\al _i)\in \al
    _i+\ndZ ^J$ for all $i\in I\setminus J$. Hence
    $\s _{i_1}\cdots \s _{i_k}^a(\al _i)\in \ndN _0^I$ for all $i\in
    I\setminus J$ by Axioms (R1) and
    (R3). Then $\ell (\s _{i_1}\cdots \s _{i_k}^a)=0$ by Lemma~\ref{le:length}
    and hence
    $\s _{i_1}\cdots \s _{i_k}^a=\id^a$.
  \end{proof}

  \begin{proposition} \label{pr:Wgfun}
    For all $J\subseteq I$ there is a unique functor
    $\Wgfun _J:\Wg (\cC |_J)\to \Wg (\cC )$
    with $\Wgfun _J(a)=a$ and $\Wgfun _J(\s _j^a)=\s _j^a $
    for all $a\in A$ and $j\in J$. This functor induces an isomorphism of
    groupoids between $\Wg (\cC |_J)$ and $\Wg _J(\cC )$.
  \end{proposition}

  \begin{proof}
    The uniqueness of $\Wgfun _J$ follows from the definition of $\Wg (\cC
    |_J)$, and $\Wgfun _J(w)\in \Wg _J(\cC )$ for all $w\in \Wg (\cC |_J)$.
    The functor $\Wgfun _J$ is well-defined by Lemma~\ref{le:wellness}.
    It is clear that $\Wgfun _J(w)=\id^a$ for some $a\in A$ and $w\in \Wg (\cC
    |_J)$ implies that $w=\id^a$, and hence $\Wgfun $ is an isomorphism.
  \end{proof}

  Finally, we state an analogue of a well-known decomposition theorem
  for Coxeter groups. Following \cite[Def.\,2.4.2]{b-BjoeBre05} let
  \begin{align}
    \Wg ^J(\cC )=\{w\in \Wg (\cC )\,|\,\ell (w \s _j)=\ell (w)+1\quad
    \text{for all $j\in J$}\}.
  \end{align}
  
  \begin{proposition} \label{pr:decomp}
    Let $J\subseteq I$ and $w\in \Wg (\cC )$. Then the following hold.
    \begin{enumerate}
      \item There exist unique elements $u \in \Wg ^J(\cC )$ and
      $v\in \Wg _J(\cC )$ such that $w=uv$.
      \item Let $u,v$ be as in (1). Then $\ell (w)=\ell (u)+\ell (v)$.
    \end{enumerate}
  \end{proposition}
  
  \begin{proof}
    The existence in (1) and the claim in (2)
    can be shown inductively on the length of $w$,
    see for example \cite[Prop.\,2.4.4]{b-BjoeBre05}. If $w\in \Wg ^J(\cC )$,
    then $w=w\id $ is a desired decomposition. Otherwise let $j\in J$ such that
    $\ell (w\s _j)=\ell (w)-1$. By induction hypothesis there exist
    $u\in \Wg ^J(\cC )$ and $v_1\in \Wg _J(\cC )$ such that $w\s _j=u v_1$
    and $\ell (w\s _j)=\ell (u)+\ell (v_1)$.
    We obtain that $w=uv$, where $v=v_1\s _j\in \Wg _J(\cC )$. Moreover
    \begin{align*}
      \ell (u)+\ell (v)&\le \ell (u)+\ell (v_1)+1=\ell (u v_1)+1\\
      &=\ell (w\s _j)+1=\ell (w)=\ell (uv)\le \ell (u)+\ell (v)
    \end{align*}
    and hence (2) holds.
    
    Let now $u_1,u_2\in \Wg ^J(\cC )$ and $v_1,v_2\in \Wg _J(\cC )$
    such that $w=u_1v_1=u_2v_2$. Then
    \begin{align} \label{eq:w1} u_1=u_2 v_2 (v_1)^{-1}. \end{align}
    Assume that $v_2\not=v_1$. Then there exists $j\in J$ such that
    $\ell (v_2 v_1^{-1}\s _j)=\ell (v_2v_1^{-1})-1$, and hence
    $v_2v_1^{-1}(\al _j)\in -\sum _{k\in J}\ndN _0\al _k$
    by Lemma~\ref{le:posi}.
    Since $u_2\in \Wg ^J(\cC )$, it follows again by Lemma~\ref{le:posi}
    that $u_2 v_2 v_1 ^{-1}(\al _j)\in -\ndN _0^I$.
    On the other hand, $u_1(\al _j)\in \ndN _0^I$
    by Lemma~\ref{le:posi} since $u_1\in \Wg ^J(\cC )$. This is
    a contradiction to \eqref{eq:w1}, and hence $v_1=v_2$ and $u_1=u_2$.
  \end{proof}
  
  An immediate consequence of Proposition~\ref{pr:decomp} is the following.
  
  \begin{corollary} \label{co:minlenrep}
    Let $J\subseteq I$. Then every left coset $w\Wg _J(\cC )$, where
    $w\in \Wg (\cC )$, has a unique representative of minimal length.
    The system of such representatives is $\Wg ^J(\cC )$.
  \end{corollary}
  
  \subsection{Geometric Combinatorics}
    Let $P$ be a partially ordered set with order relation $\preceq$. 
    A \textit{chain} of length $i$ in $P$ is a linearly ordered subset  
    $p_0 \prec \cdots \prec p_i$ of $i+1$ elements of $P$. 
    A chain is called \textit{maximal} if it is an inclusionwise maximal linearly
    ordered subset of $P$.  
    The order complex $\Delta(P)$ of $P$ is the
    abstract simplicial complex on ground set $P$ whose $i$-simplices
    are the chains of length $i$.  If $p \preceq q$ are two elements of
    $P$ then we denote by
    $[p,q]$ the {\it closed interval} $\{ r \in P~|~p \preceq r \preceq q \}$.
    Analogously, one defines the
    {\it open interval} $(p,q) := [p,q] \setminus \{p,q\}$.
    We write $\Delta(p,q)$ to denote the {\it order complex} of $(p,q)$.
    For $p \in P$ we write $P_{\prec p}$ for the subposet of all 
    $q \in P$ with $q \prec p$. 
  
    Via the geometric realization $|\Delta(P)|$ of $P$ one can speak of
    topological properties of
    partially ordered sets $P$. In particular, we can speak of $P$ being
    homotopy equivalent or
    homeomorphic to another partially ordered set or topological space.
    If $P$ is a partially ordered set with unique maximal element $\hat{1}$ or
    unique minimal element $\hat{0}$ then $\Delta(P)$ is a cone over
     $\hat{1}$ (resp. $\hat{0}$) and therefore contractible. Hence in order to
    be able to capture non-trivial topology one considers for partially ordered
    sets $P$ with unique minimal element $\hat{0}$ and unique maximal element 
    $\hat{1}$ the proper part $\hat{P} := P \setminus \{ \hat{0}, \hat{1}\}$
    of $P$. For example $\widehat{[p,q]} = (p,q)$. The following simple
    example will be useful in the subsequent sections.

    \begin{example} \label{ex:boolean}
      Let $\Omega$ be a finite set and $2^\Omega$ be the Boolean lattice of
      all subsets of
      $\Omega$ ordered by inclusion. Then $2^\Omega$ has
      unique minimal element $\hat{0} = \emptyset$ and unique maximal element
      $\hat{1} = \Omega$. Then $\Delta(\widehat{2^\Omega})$ is the 
      barycentric subdivision (see for example \cite[\S 15]{b-Munkres84}) 
      of the boundary of the $(|\Omega|-1)$-simplex and hence
      homeomorphic to an $(|\Omega|-2)$-sphere.
    \end{example}

    For our purposes the following well known result on the topology of
    partially ordered sets will be crucial.

    \begin{theorem}[Corollary 10.12 \cite{a-Bj95}] \label{th:closure}
       Let $P$ be a partially ordered set and let $f : P \rightarrow P$ be
       a map such that:
       \begin{itemize}
          \item[(1)] $p \preceq q$ implies $f(p) \preceq f(q)$.
          \item[(2)] $f(p) \preceq p$.
       \end{itemize}
       Then $P$ and $f(P)$ are homotopy equivalent.
    \end{theorem}

    In order to set up the next tool it is most convenient to work in the context of (abstract)
    simplicial complexes. For a simplicial complex $\Delta$ we call $A \in 
    \Delta$ a \textit{face} of $\Delta$ and denote by $\dim A = \# A -1$ its \textit{dimension}. 
    We call $\Delta$ \textit{pure} if all inclusionwise maximal faces have the same dimension. 
    The order complex $\Delta(P)$ of a partially ordered set $P$ is pure if and
    only if all maximal chains in $P$ have the same length. 
    A pure simplicial complex $\Delta$ is called \textit{shellable} if there is a
    numbering $F_1,  \ldots , F_r$ of the set of its maximal faces such that
    for all $1 \leq i < j \leq r$ there is an $\ell < j$ and an $\omega \in
    F_j$ such that $F_i \cap F_j \subseteq  F_\ell \cap F_j = F_j \setminus
    \{ \omega\}$. 

    It is well known (see e.g. \cite{a-Bj95}) that if $\Delta$ is a shellable simplicial complex of
    dimension $d$ then the geometric realization is homotopy equivalent to a
    wedge of spheres of dimension $d$. For the subsequent applications we are
    interested in situations when $\Delta$ is homeomorphic to a sphere. 
    This can also be verified using shellability when $\Delta$ is a
    pseudmanifold. A pure $d$-dimensional simplicial complex $\Delta$ is 
    called a \textit{pseudomanifold} if for all faces $F \in \Delta$ of dimension $d-1$
    there are at most $2$ faces of dimension $d$ containing $F$.

    \begin{theorem}[Theorem 11.4 \cite{a-Bj95}]
      \label{th:plsphere}
      Let $\Delta$ be a shellable $d$-dimensional pseudomanifold. If every face of 
      dimension $d-1$ is contained in exactly $2$ faces of dimension $d$ then
      $\Delta$ is homeomorphic to a $d$-sphere otherwise $\Delta$ is
      homeomorphic to a $d$-ball.
   \end{theorem}

  \section{Weak Order}
  \label{sec:weak}

  Throughout this section let
  $\cC =\cC (I,A,(\rfl _i)_{i\in I},(\Cm ^a)_{a\in A})$ be a
  Cartan scheme
  and assume that $\rsC \re (\cC )$ is a finite root system.

  The (right)
  \textit{weak order} or \textit{Duflo order}
  on Weyl groupoids is the natural generalization of the (right)
  weak order on Coxeter
  groups, see \cite[Ch.\,3]{b-BjoeBre05}:
  for any $a,b,c\in A$ and $u\in \Hom(b,a)$,
  $v\in \Hom (c,b)$ we define
  $$u\wo uv~: \Leftrightarrow~\ell (u)+\ell (v) = \ell (uv).$$
  For all $a\in A$ the weak order is a partial ordering on $\Homsto a$.
  As shown in \cite{p-HeckSchn09}, the weak order
  has an algebraic interpretation
  in terms of right coideal subalgebras of Nichols algebras.

  \begin{example} \label{ex:bruhat}
    Let $I=\{1,2,3\}$ and $A=\{a,b,c,d,e\}$. There is a unique Cartan
    scheme $\cC $ with
    \begin{gather*}
      C^a=
      \begin{pmatrix}
        2 & -1 & 0 \\
        -1 & 2 & -2 \\
        0 & -1 & 2
      \end{pmatrix}, \quad
      C^b=
      \begin{pmatrix}
        2 & -1 & 0 \\
        -1 & 2 & -1 \\
        0 & -1 & 2
      \end{pmatrix}, \quad
      C^c=
      \begin{pmatrix}
        2 & -1 & -1 \\
        -1 & 2 & -1 \\
        -1 & -1 & 2
      \end{pmatrix},\\
      C^d=
      \begin{pmatrix}
        2 & 0 & -1 \\
        0 & 2 & -1 \\
        -1 & -1 & 2
      \end{pmatrix}, \quad
      C^e=
      \begin{pmatrix}
        2 & 0 & -1 \\
        0 & 2 & -1 \\
        -1 & -2 & 2
      \end{pmatrix},
    \end{gather*}
    where the object change diagram is as in Figure~\ref{fi:ocd8}.
    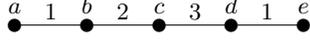
\begin{figure}
      \setlength{\unitlength}{.8mm}
      \begin{picture}(50,8)
        \put(1,2){\circle*{2}}
        \put(1,4){\makebox[0pt]{\scriptsize $a$}}
        \put(2,2){\line(1,0){10}}
        \put(13,2){\circle*{2}}
        \put(13,4){\makebox[0pt]{\scriptsize $b$}}
        \put(14,2){\line(1,0){10}}
        \put(25,2){\circle*{2}}
        \put(25,4){\makebox[0pt]{\scriptsize $c$}}
        \put(26,2){\line(1,0){10}}
        \put(37,2){\circle*{2}}
        \put(37,4){\makebox[0pt]{\scriptsize $d$}}
        \put(38,2){\line(1,0){10}}
        \put(49,2){\circle*{2}}
        \put(49,4){\makebox[0pt]{\scriptsize $e$}}
        \put(7,3){\makebox[0pt]{\scriptsize $1$}}
        \put(19,3){\makebox[0pt]{\scriptsize $2$}}
        \put(31,3){\makebox[0pt]{\scriptsize $3$}}
        \put(43,3){\makebox[0pt]{\scriptsize $1$}}
      \end{picture}
      \caption{The object change diagram for Example \ref{ex:bruhat}}
      \label{fi:ocd8}
    \end{figure}

    \begin{figure}

      \setlength{\unitlength}{1973sp}%

      \vskip11cm
      \hskip-8cm\begin{picture}(0,0)
       \includegraphics{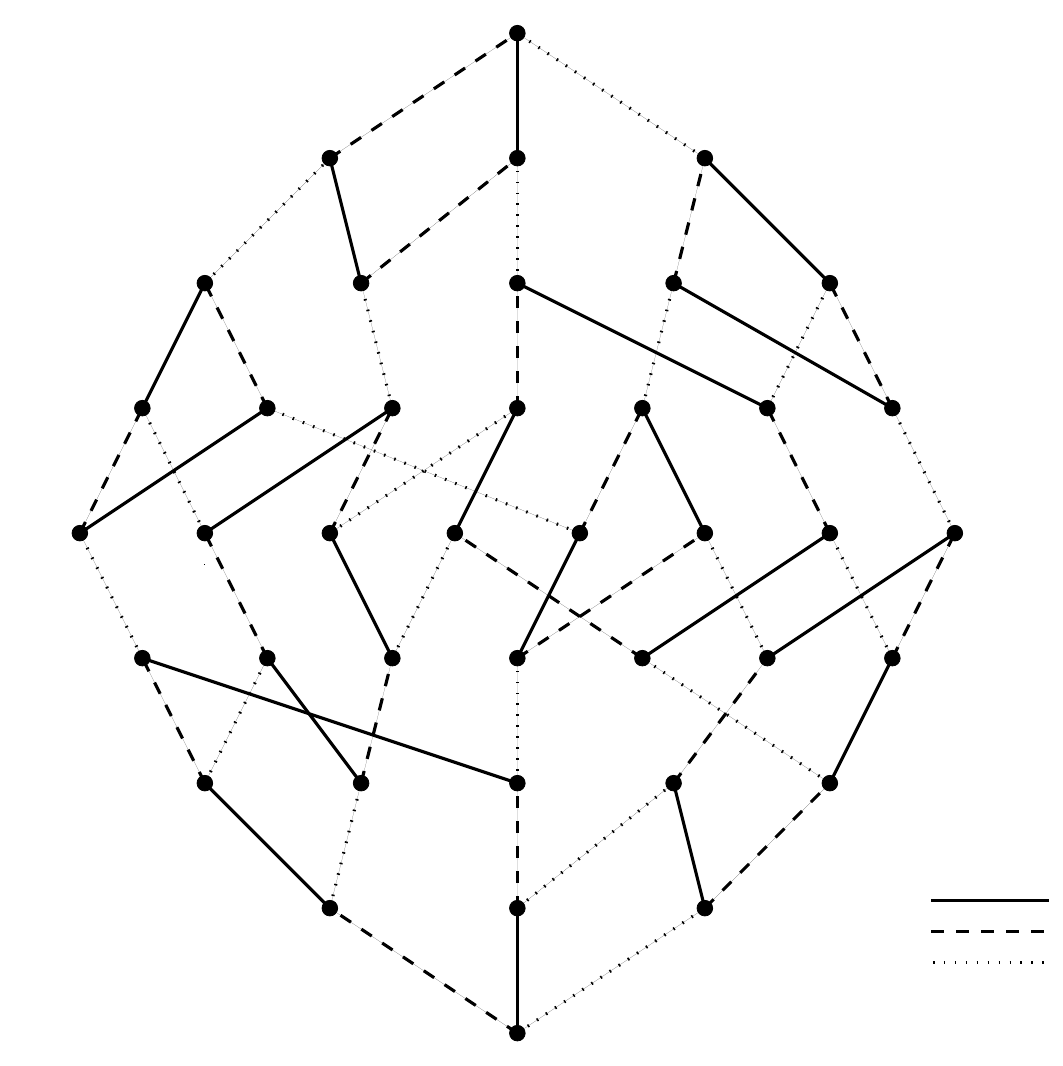}%
      \end{picture}%

\everymath{\scriptstyle}
      \begingroup\makeatletter\ifx\SetFigFont\undefined%
         \gdef\SetFigFont#1#2#3#4#5{%
         \reset@font\fontsize{#1}{#2pt}%
         \fontfamily{#3}\fontseries{#4}\fontshape{#5}%
         \selectfont}%
      \fi\endgroup%
      \hskip8cm\vskip-11cm
      \begin{picture}(10068,10230)(736,-9601)
         \put(10600,-7700){\makebox(0,0)[lb]{\smash{{\SetFigFont{12}{14.4}{\familydefault}{\mddefault}{\updefault}{\color[rgb]{0,0,0}${1}$}%
}}}}
         \put(10600,-8100){\makebox(0,0)[lb]{\smash{{\SetFigFont{12}{14.4}{\familydefault}{\mddefault}{\updefault}{\color[rgb]{0,0,0}${2}$}%
}}}}
         \put(10600,-8500){\makebox(0,0)[lb]{\smash{{\SetFigFont{12}{14.4}{\familydefault}{\mddefault}{\updefault}{\color[rgb]{0,0,0}${3}$}%
}}}}

         \put(6900,-9400){\makebox(0,0)[lb]{\smash{{\SetFigFont{12}{14.4}{\familydefault}{\mddefault}{\updefault}{\color[rgb]{0,0,0}${\id^a}$}%
}}}}

         \put(4900,-8050){\makebox(0,0)[lb]{\smash{{\SetFigFont{12}{14.4}{\familydefault}{\mddefault}{\updefault}{\color[rgb]{0,0,0}${2^a}$}%
}}}}
         \put(7000,-8100){\makebox(0,0)[lb]{\smash{{\SetFigFont{12}{14.4}{\familydefault}{\mddefault}{\updefault}{\color[rgb]{0,0,0}${1^b}$}%
}}}}
         \put(8800,-8100){\makebox(0,0)[lb]{\smash{{\SetFigFont{12}{14.4}{\familydefault}{\mddefault}{\updefault}{\color[rgb]{0,0,0}${3^a}$}%
}}}}

         \put(3500,-6950){\makebox(0,0)[lb]{\smash{{\SetFigFont{12}{14.4}{\familydefault}{\mddefault}{\updefault}{\color[rgb]{0,0,0}${21^b}$}%
}}}}
         \put(4900,-6850){\makebox(0,0)[lb]{\smash{{\SetFigFont{12}{14.4}{\familydefault}{\mddefault}{\updefault}{\color[rgb]{0,0,0}${23^a}$}%
}}}}
         \put(7050,-6850){\makebox(0,0)[lb]{\smash{{\SetFigFont{12}{14.4}{\familydefault}{\mddefault}{\updefault}{\color[rgb]{0,0,0}${12^c}$}%
}}}}
         \put(8600,-6750){\makebox(0,0)[lb]{\smash{{\SetFigFont{12}{14.4}{\familydefault}{\mddefault}{\updefault}{\color[rgb]{0,0,0}${13^b}$}%
}}}}
         \put(9800,-6950){\makebox(0,0)[lb]{\smash{{\SetFigFont{12}{14.4}{\familydefault}{\mddefault}{\updefault}{\color[rgb]{0,0,0}${32^a}$}%
}}}}

         \put(2600,-5700){\makebox(0,0)[lb]{\smash{{\SetFigFont{12}{14.4}{\familydefault}{\mddefault}{\updefault}{\color[rgb]{0,0,0}${121^c}$}%
}}}}
         \put(3800,-5300){\makebox(0,0)[lb]{\smash{{\SetFigFont{12}{14.4}{\familydefault}{\mddefault}{\updefault}{\color[rgb]{0,0,0}${213^b}$}%
}}}}
         \put(5000,-5650){\makebox(0,0)[lb]{\smash{{\SetFigFont{12}{14.4}{\familydefault}{\mddefault}{\updefault}{\color[rgb]{0,0,0}${232^a}$}%
}}}}
         \put(6200,-5700){\makebox(0,0)[lb]{\smash{{\SetFigFont{12}{14.4}{\familydefault}{\mddefault}{\updefault}{\color[rgb]{0,0,0}${123^d}$}%
}}}}
         \put(7600,-5700){\makebox(0,0)[lb]{\smash{{\SetFigFont{12}{14.4}{\familydefault}{\mddefault}{\updefault}{\color[rgb]{0,0,0}${323^a}$}%
}}}}
         \put(9400,-5700){\makebox(0,0)[lb]{\smash{{\SetFigFont{12}{14.4}{\familydefault}{\mddefault}{\updefault}{\color[rgb]{0,0,0}${132^c}$}%
}}}}
         \put(10600,-5700){\makebox(0,0)[lb]{\smash{{\SetFigFont{12}{14.4}{\familydefault}{\mddefault}{\updefault}{\color[rgb]{0,0,0}${321^b}$}%
}}}}

         \put(1900,-4400){\makebox(0,0)[lb]{\smash{{\SetFigFont{12}{14.4}{\familydefault}{\mddefault}{\updefault}{\color[rgb]{0,0,0}${1213^d}$}%
}}}}
         \put(3200,-4400){\makebox(0,0)[lb]{\smash{{\SetFigFont{12}{14.4}{\familydefault}{\mddefault}{\updefault}{\color[rgb]{0,0,0}${2132^c}$}%
}}}}
         \put(4350,-4400){\makebox(0,0)[lb]{\smash{{\SetFigFont{12}{14.4}{\familydefault}{\mddefault}{\updefault}{\color[rgb]{0,0,0}${2321^b}$}%
}}}}
         \put(5500,-4400){\makebox(0,0)[lb]{\smash{{\SetFigFont{12}{14.4}{\familydefault}{\mddefault}{\updefault}{\color[rgb]{0,0,0}${2323^a}$}%
}}}}
         \put(6750,-4400){\makebox(0,0)[lb]{\smash{{\SetFigFont{12}{14.4}{\familydefault}{\mddefault}{\updefault}{\color[rgb]{0,0,0}${1231^e}$}%
}}}}
         \put(7950,-4200){\makebox(0,0)[lb]{\smash{{\SetFigFont{12}{14.4}{\familydefault}{\mddefault}{\updefault}{\color[rgb]{0,0,0}${1232^d}$}%
}}}}
         \put(10050,-4300){\makebox(0,0)[lb]{\smash{{\SetFigFont{12}{14.4}{\familydefault}{\mddefault}{\updefault}{\color[rgb]{0,0,0}${3213^b}$}%
}}}}
         \put(11200,-4400){\makebox(0,0)[lb]{\smash{{\SetFigFont{12}{14.4}{\familydefault}{\mddefault}{\updefault}{\color[rgb]{0,0,0}${1321^c}$}%
}}}}

         \put(2350,-3200){\makebox(0,0)[lb]{\smash{{\SetFigFont{12}{14.4}{\familydefault}{\mddefault}{\updefault}{\color[rgb]{0,0,0}$12132^d$}%
}}}}
         \put(3550,-3050){\makebox(0,0)[lb]{\smash{{\SetFigFont{12}{14.4}{\familydefault}{\mddefault}{\updefault}{\color[rgb]{0,0,0}$12131^e$}%
}}}}
         \put(4800,-3000){\makebox(0,0)[lb]{\smash{{\SetFigFont{12}{14.4}{\familydefault}{\mddefault}{\updefault}{\color[rgb]{0,0,0}$21321^c$}%
}}}}
         \put(6000,-3000){\makebox(0,0)[lb]{\smash{{\SetFigFont{12}{14.4}{\familydefault}{\mddefault}{\updefault}{\color[rgb]{0,0,0}$23213^b$}%
}}}}
         \put(7200,-3100){\makebox(0,0)[lb]{\smash{{\SetFigFont{12}{14.4}{\familydefault}{\mddefault}{\updefault}{\color[rgb]{0,0,0}$12312^e$}%
}}}}
         \put(8500,-3200){\makebox(0,0)[lb]{\smash{{\SetFigFont{12}{14.4}{\familydefault}{\mddefault}{\updefault}{\color[rgb]{0,0,0}$32132^c$}%
}}}}
         \put(10700,-3200){\makebox(0,0)[lb]{\smash{{\SetFigFont{12}{14.4}{\familydefault}{\mddefault}{\updefault}{\color[rgb]{0,0,0}$13213^d$}%
}}}}

         \put(2850,-1900){\makebox(0,0)[lb]{\smash{{\SetFigFont{12}{14.4}{\familydefault}{\mddefault}{\updefault}{\color[rgb]{0,0,0}$121312^e$}%
}}}}
         \put(4350,-1900){\makebox(0,0)[lb]{\smash{{\SetFigFont{12}{14.4}{\familydefault}{\mddefault}{\updefault}{\color[rgb]{0,0,0}$213213^d$}%
}}}}
         \put(5850,-1900){\makebox(0,0)[lb]{\smash{{\SetFigFont{12}{14.4}{\familydefault}{\mddefault}{\updefault}{\color[rgb]{0,0,0}$232132^c$}%
}}}}
         \put(7350,-1900){\makebox(0,0)[lb]{\smash{{\SetFigFont{12}{14.4}{\familydefault}{\mddefault}{\updefault}{\color[rgb]{0,0,0}$123123^e$}%
}}}}
         \put(10100,-1900){\makebox(0,0)[lb]{\smash{{\SetFigFont{12}{14.4}{\familydefault}{\mddefault}{\updefault}{\color[rgb]{0,0,0}$132132^d$}%
}}}}

         \put(3850,-600){\makebox(0,0)[lb]{\smash{{\SetFigFont{12}{14.4}{\familydefault}{\mddefault}{\updefault}{\color[rgb]{0,0,0}$1213123^e$}%
}}}}
         \put(5650,-600){\makebox(0,0)[lb]{\smash{{\SetFigFont{12}{14.4}{\familydefault}{\mddefault}{\updefault}{\color[rgb]{0,0,0}$2132132^d$}%
}}}}
         \put(8830,-600){\makebox(0,0)[lb]{\smash{{\SetFigFont{12}{14.4}{\familydefault}{\mddefault}{\updefault}{\color[rgb]{0,0,0}$1231232^e$}%
}}}}

         \put(6400,800){\makebox(0,0)[lb]{\smash{{\SetFigFont{12}{14.4}{\familydefault}{\mddefault}{\updefault}{\color[rgb]{0,0,0}${12131232^e}$}%
}}}}
       \end{picture}%
       \caption{The weak order for Example \ref{ex:bruhat} in object a}
       \label{fi:bruhata}
    \end{figure}

    \begin{figure}

      \setlength{\unitlength}{1973sp}%

      \vskip11cm
      \hskip-6.5cm
      \begin{picture}(0,0)
       \includegraphics{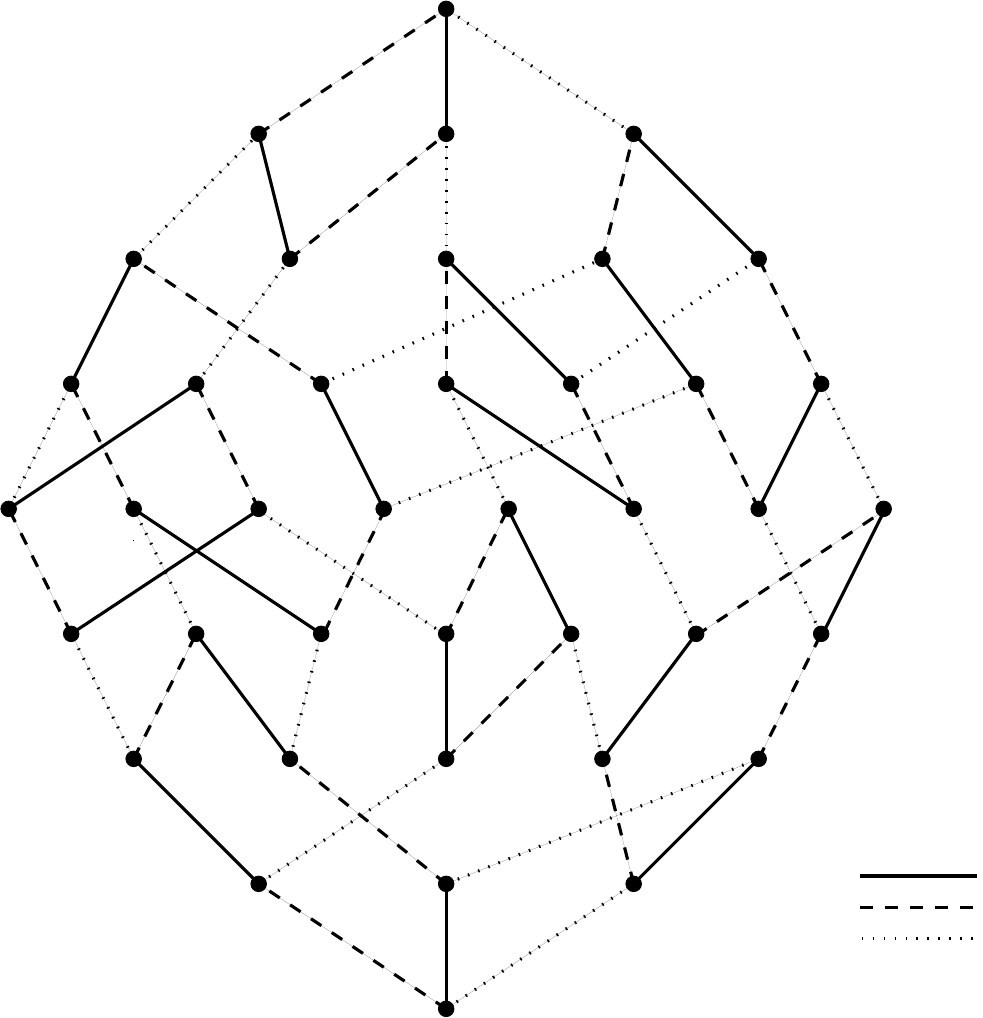}%
      \end{picture}%

      \everymath{\scriptstyle}
      \begingroup\makeatletter\ifx\SetFigFont\undefined%
         \gdef\SetFigFont#1#2#3#4#5{%
         \reset@font\fontsize{#1}{#2pt}%
         \fontfamily{#3}\fontseries{#4}\fontshape{#5}%
         \selectfont}%
      \fi\endgroup%
      \hskip8cm\vskip-10.8cm
      \begin{picture}(10068,10230)(736,-9601)

         \put(10600,-7700){\makebox(0,0)[lb]{\smash{{\SetFigFont{12}{14.4}{\familydefault}{\mddefault}{\updefault}{\color[rgb]{0,0,0}${1}$}%
}}}}
         \put(10600,-8100){\makebox(0,0)[lb]{\smash{{\SetFigFont{12}{14.4}{\familydefault}{\mddefault}{\updefault}{\color[rgb]{0,0,0}${2}$}%
}}}}
         \put(10600,-8500){\makebox(0,0)[lb]{\smash{{\SetFigFont{12}{14.4}{\familydefault}{\mddefault}{\updefault}{\color[rgb]{0,0,0}${3}$}%
}}}}

         \put(6900,-9400){\makebox(0,0)[lb]{\smash{{\SetFigFont{12}{14.4}{\familydefault}{\mddefault}{\updefault}{\color[rgb]{0,0,0}${\id^b}$}%
}}}}

         \put(4900,-8050){\makebox(0,0)[lb]{\smash{{\SetFigFont{12}{14.4}{\familydefault}{\mddefault}{\updefault}{\color[rgb]{0,0,0}$2^c$}%
}}}}
         \put(7000,-8100){\makebox(0,0)[lb]{\smash{{\SetFigFont{12}{14.4}{\familydefault}{\mddefault}{\updefault}{\color[rgb]{0,0,0}$1^a$}%
}}}}
         \put(8800,-8100){\makebox(0,0)[lb]{\smash{{\SetFigFont{12}{14.4}{\familydefault}{\mddefault}{\updefault}{\color[rgb]{0,0,0}$3^b$}%
}}}}

         \put(3500,-6950){\makebox(0,0)[lb]{\smash{{\SetFigFont{12}{14.4}{\familydefault}{\mddefault}{\updefault}{\color[rgb]{0,0,0}$21^c$}%
}}}}
         \put(4900,-6850){\makebox(0,0)[lb]{\smash{{\SetFigFont{12}{14.4}{\familydefault}{\mddefault}{\updefault}{\color[rgb]{0,0,0}$12^a$}%
}}}}
         \put(7050,-6850){\makebox(0,0)[lb]{\smash{{\SetFigFont{12}{14.4}{\familydefault}{\mddefault}{\updefault}{\color[rgb]{0,0,0}$23^d$}%
}}}}
         \put(8600,-6750){\makebox(0,0)[lb]{\smash{{\SetFigFont{12}{14.4}{\familydefault}{\mddefault}{\updefault}{\color[rgb]{0,0,0}$32^c$}%
}}}}
         \put(9800,-6950){\makebox(0,0)[lb]{\smash{{\SetFigFont{12}{14.4}{\familydefault}{\mddefault}{\updefault}{\color[rgb]{0,0,0}$13^a$}%
}}}}

         \put(2700,-5700){\makebox(0,0)[lb]{\smash{{\SetFigFont{12}{14.4}{\familydefault}{\mddefault}{\updefault}{\color[rgb]{0,0,0}$213^d$}%
}}}}
         \put(3950,-5500){\makebox(0,0)[lb]{\smash{{\SetFigFont{12}{14.4}{\familydefault}{\mddefault}{\updefault}{\color[rgb]{0,0,0}$121^b$}%
}}}}
         \put(5100,-5650){\makebox(0,0)[lb]{\smash{{\SetFigFont{12}{14.4}{\familydefault}{\mddefault}{\updefault}{\color[rgb]{0,0,0}$123^a$}%
}}}}
         \put(6350,-5700){\makebox(0,0)[lb]{\smash{{\SetFigFont{12}{14.4}{\familydefault}{\mddefault}{\updefault}{\color[rgb]{0,0,0}$231^e$}%
}}}}
         \put(8300,-5600){\makebox(0,0)[lb]{\smash{{\SetFigFont{12}{14.4}{\familydefault}{\mddefault}{\updefault}{\color[rgb]{0,0,0}$232^d$}%
}}}}
         \put(9400,-5700){\makebox(0,0)[lb]{\smash{{\SetFigFont{12}{14.4}{\familydefault}{\mddefault}{\updefault}{\color[rgb]{0,0,0}$321^c$}%
}}}}
         \put(10600,-5700){\makebox(0,0)[lb]{\smash{{\SetFigFont{12}{14.4}{\familydefault}{\mddefault}{\updefault}{\color[rgb]{0,0,0}$132^a$}%
}}}}

         \put(1900,-4400){\makebox(0,0)[lb]{\smash{{\SetFigFont{12}{14.4}{\familydefault}{\mddefault}{\updefault}{\color[rgb]{0,0,0}$2132^d$}%
}}}}
         \put(3200,-4400){\makebox(0,0)[lb]{\smash{{\SetFigFont{12}{14.4}{\familydefault}{\mddefault}{\updefault}{\color[rgb]{0,0,0}$1213^b$}%
}}}}
         \put(4350,-4300){\makebox(0,0)[lb]{\smash{{\SetFigFont{12}{14.4}{\familydefault}{\mddefault}{\updefault}{\color[rgb]{0,0,0}$2131^e$}%
}}}}
         \put(5550,-4400){\makebox(0,0)[lb]{\smash{{\SetFigFont{12}{14.4}{\familydefault}{\mddefault}{\updefault}{\color[rgb]{0,0,0}$1232^a$}%
}}}}
         \put(6750,-4400){\makebox(0,0)[lb]{\smash{{\SetFigFont{12}{14.4}{\familydefault}{\mddefault}{\updefault}{\color[rgb]{0,0,0}$2312^e$}%
}}}}
         \put(8000,-4450){\makebox(0,0)[lb]{\smash{{\SetFigFont{12}{14.4}{\familydefault}{\mddefault}{\updefault}{\color[rgb]{0,0,0}$3213^d$}%
}}}}
         \put(10050,-4300){\makebox(0,0)[lb]{\smash{{\SetFigFont{12}{14.4}{\familydefault}{\mddefault}{\updefault}{\color[rgb]{0,0,0}$1323^a$}%
}}}}
         \put(11250,-4400){\makebox(0,0)[lb]{\smash{{\SetFigFont{12}{14.4}{\familydefault}{\mddefault}{\updefault}{\color[rgb]{0,0,0}$1321^b$}%
}}}}

         \put(2350,-3200){\makebox(0,0)[lb]{\smash{{\SetFigFont{12}{14.4}{\familydefault}{\mddefault}{\updefault}{\color[rgb]{0,0,0}$12132^c$}%
}}}}
         \put(3600,-3050){\makebox(0,0)[lb]{\smash{{\SetFigFont{12}{14.4}{\familydefault}{\mddefault}{\updefault}{\color[rgb]{0,0,0}$21312^e$}%
}}}}
         \put(4900,-3250){\makebox(0,0)[lb]{\smash{{\SetFigFont{12}{14.4}{\familydefault}{\mddefault}{\updefault}{\color[rgb]{0,0,0}$12321^b$}%
}}}}
         \put(6070,-3100){\makebox(0,0)[lb]{\smash{{\SetFigFont{12}{14.4}{\familydefault}{\mddefault}{\updefault}{\color[rgb]{0,0,0}$23123^e$}%
}}}}
         \put(7240,-3150){\makebox(0,0)[lb]{\smash{{\SetFigFont{12}{14.4}{\familydefault}{\mddefault}{\updefault}{\color[rgb]{0,0,0}$32132^d$}%
}}}}
         \put(9400,-2900){\makebox(0,0)[lb]{\smash{{\SetFigFont{12}{14.4}{\familydefault}{\mddefault}{\updefault}{\color[rgb]{0,0,0}$12323^a$}%
}}}}
         \put(10700,-3200){\makebox(0,0)[lb]{\smash{{\SetFigFont{12}{14.4}{\familydefault}{\mddefault}{\updefault}{\color[rgb]{0,0,0}$13213^b$}%
}}}}

         \put(2850,-1900){\makebox(0,0)[lb]{\smash{{\SetFigFont{12}{14.4}{\familydefault}{\mddefault}{\updefault}{\color[rgb]{0,0,0}$121321^c$}%
}}}}
         \put(4350,-1900){\makebox(0,0)[lb]{\smash{{\SetFigFont{12}{14.4}{\familydefault}{\mddefault}{\updefault}{\color[rgb]{0,0,0}$213123^e$}%
}}}}
         \put(5850,-1900){\makebox(0,0)[lb]{\smash{{\SetFigFont{12}{14.4}{\familydefault}{\mddefault}{\updefault}{\color[rgb]{0,0,0}$231232^e$}%
}}}}
         \put(7350,-1900){\makebox(0,0)[lb]{\smash{{\SetFigFont{12}{14.4}{\familydefault}{\mddefault}{\updefault}{\color[rgb]{0,0,0}$123213^b$}%
}}}}
         \put(10100,-1900){\makebox(0,0)[lb]{\smash{{\SetFigFont{12}{14.4}{\familydefault}{\mddefault}{\updefault}{\color[rgb]{0,0,0}$132132^c$}%
}}}}

         \put(3870,-600){\makebox(0,0)[lb]{\smash{{\SetFigFont{12}{14.4}{\familydefault}{\mddefault}{\updefault}{\color[rgb]{0,0,0}$1213213^d$}%
}}}}
         \put(5720,-600){\makebox(0,0)[lb]{\smash{{\SetFigFont{12}{14.4}{\familydefault}{\mddefault}{\updefault}{\color[rgb]{0,0,0}$2131232^e$}%
}}}}
         \put(8750,-500){\makebox(0,0)[lb]{\smash{{\SetFigFont{12}{14.4}{\familydefault}{\mddefault}{\updefault}{\color[rgb]{0,0,0}$1232132^c$}%
}}}}

         \put(6400,800){\makebox(0,0)[lb]{\smash{{\SetFigFont{12}{14.4}{\familydefault}{\mddefault}{\updefault}{\color[rgb]{0,0,0}$12132132^d$}%
}}}}
       \end{picture}%
       \caption{The weak order for Example \ref{ex:bruhat} in object b}
       \label{fi:bruhatb}
    \end{figure}

    \begin{figure}

      \setlength{\unitlength}{1973sp}%

      \vskip11cm
      \hskip-6.5cm
      \begin{picture}(0,0)
       \includegraphics{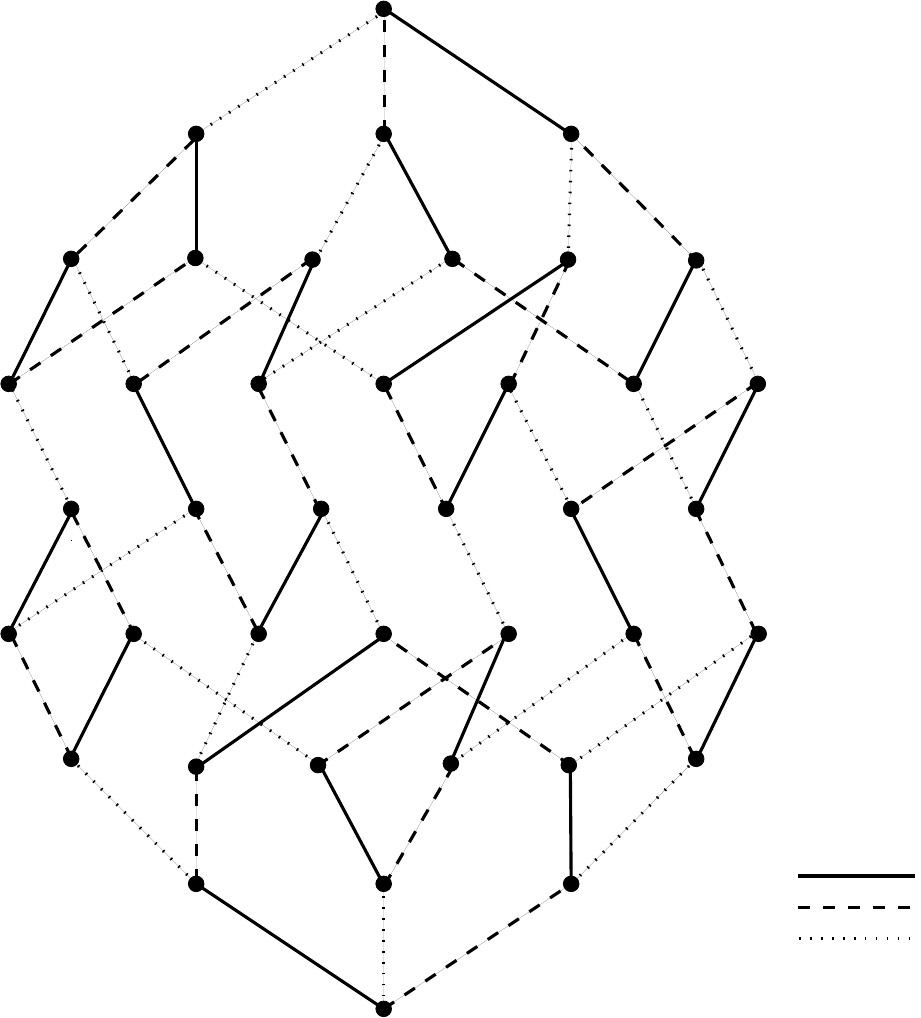}%
      \end{picture}%

      \everymath{\scriptstyle}
      \begingroup\makeatletter\ifx\SetFigFont\undefined%
         \gdef\SetFigFont#1#2#3#4#5{%
         \reset@font\fontsize{#1}{#2pt}%
         \fontfamily{#3}\fontseries{#4}\fontshape{#5}%
         \selectfont}%
      \fi\endgroup%
      \hskip8cm\vskip-10.8cm
      \begin{picture}(10068,10230)(736,-9601)

         \put(10100,-7700){\makebox(0,0)[lb]{\smash{{\SetFigFont{12}{14.4}{\familydefault}{\mddefault}{\updefault}{\color[rgb]{0,0,0}${1}$}%
}}}}
         \put(10100,-8100){\makebox(0,0)[lb]{\smash{{\SetFigFont{12}{14.4}{\familydefault}{\mddefault}{\updefault}{\color[rgb]{0,0,0}${2}$}%
}}}}
         \put(10100,-8500){\makebox(0,0)[lb]{\smash{{\SetFigFont{12}{14.4}{\familydefault}{\mddefault}{\updefault}{\color[rgb]{0,0,0}${3}$}%
}}}}

         \put(6300,-9400){\makebox(0,0)[lb]{\smash{{\SetFigFont{12}{14.4}{\familydefault}{\mddefault}{\updefault}{\color[rgb]{0,0,0}${\id^c}$}%
}}}}

         \put(4250,-8050){\makebox(0,0)[lb]{\smash{{\SetFigFont{12}{14.4}{\familydefault}{\mddefault}{\updefault}{\color[rgb]{0,0,0}$1^c$}%
}}}}
         \put(6500,-8100){\makebox(0,0)[lb]{\smash{{\SetFigFont{12}{14.4}{\familydefault}{\mddefault}{\updefault}{\color[rgb]{0,0,0}$3^d$}%
}}}}
         \put(8200,-8100){\makebox(0,0)[lb]{\smash{{\SetFigFont{12}{14.4}{\familydefault}{\mddefault}{\updefault}{\color[rgb]{0,0,0}$2^b$}%
}}}}

         \put(2900,-6950){\makebox(0,0)[lb]{\smash{{\SetFigFont{12}{14.4}{\familydefault}{\mddefault}{\updefault}{\color[rgb]{0,0,0}$13^d$}%
}}}}
         \put(4050,-6850){\makebox(0,0)[lb]{\smash{{\SetFigFont{12}{14.4}{\familydefault}{\mddefault}{\updefault}{\color[rgb]{0,0,0}$12^b$}%
}}}}
         \put(5250,-6850){\makebox(0,0)[lb]{\smash{{\SetFigFont{12}{14.4}{\familydefault}{\mddefault}{\updefault}{\color[rgb]{0,0,0}$31^e$}%
}}}}
         \put(6550,-6750){\makebox(0,0)[lb]{\smash{{\SetFigFont{12}{14.4}{\familydefault}{\mddefault}{\updefault}{\color[rgb]{0,0,0}$32^d$}%
}}}}
         \put(8300,-6950){\makebox(0,0)[lb]{\smash{{\SetFigFont{12}{14.4}{\familydefault}{\mddefault}{\updefault}{\color[rgb]{0,0,0}$21^a$}%
}}}}
         \put(9450,-6950){\makebox(0,0)[lb]{\smash{{\SetFigFont{12}{14.4}{\familydefault}{\mddefault}{\updefault}{\color[rgb]{0,0,0}$23^b$}%
}}}}

         \put(2100,-5700){\makebox(0,0)[lb]{\smash{{\SetFigFont{12}{14.4}{\familydefault}{\mddefault}{\updefault}{\color[rgb]{0,0,0}$132^d$}%
}}}}
         \put(3300,-5500){\makebox(0,0)[lb]{\smash{{\SetFigFont{12}{14.4}{\familydefault}{\mddefault}{\updefault}{\color[rgb]{0,0,0}$131^e$}%
}}}}
         \put(4450,-5500){\makebox(0,0)[lb]{\smash{{\SetFigFont{12}{14.4}{\familydefault}{\mddefault}{\updefault}{\color[rgb]{0,0,0}$123^b$}%
}}}}
         \put(6350,-5250){\makebox(0,0)[lb]{\smash{{\SetFigFont{12}{14.4}{\familydefault}{\mddefault}{\updefault}{\color[rgb]{0,0,0}$121^a$}%
}}}}
         \put(7600,-5250){\makebox(0,0)[lb]{\smash{{\SetFigFont{12}{14.4}{\familydefault}{\mddefault}{\updefault}{\color[rgb]{0,0,0}$312^e$}%
}}}}
         \put(8950,-5500){\makebox(0,0)[lb]{\smash{{\SetFigFont{12}{14.4}{\familydefault}{\mddefault}{\updefault}{\color[rgb]{0,0,0}$232^c$}%
}}}}
         \put(10100,-5700){\makebox(0,0)[lb]{\smash{{\SetFigFont{12}{14.4}{\familydefault}{\mddefault}{\updefault}{\color[rgb]{0,0,0}$213^a$}%
}}}}

         \put(2500,-4200){\makebox(0,0)[lb]{\smash{{\SetFigFont{12}{14.4}{\familydefault}{\mddefault}{\updefault}{\color[rgb]{0,0,0}$1312^e$}%
}}}}
         \put(3700,-4200){\makebox(0,0)[lb]{\smash{{\SetFigFont{12}{14.4}{\familydefault}{\mddefault}{\updefault}{\color[rgb]{0,0,0}$1232^c$}%
}}}}
         \put(4900,-4200){\makebox(0,0)[lb]{\smash{{\SetFigFont{12}{14.4}{\familydefault}{\mddefault}{\updefault}{\color[rgb]{0,0,0}$1213^a$}%
}}}}
         \put(6100,-4200){\makebox(0,0)[lb]{\smash{{\SetFigFont{12}{14.4}{\familydefault}{\mddefault}{\updefault}{\color[rgb]{0,0,0}$3123^e$}%
}}}}
         \put(7350,-4200){\makebox(0,0)[lb]{\smash{{\SetFigFont{12}{14.4}{\familydefault}{\mddefault}{\updefault}{\color[rgb]{0,0,0}$2321^c$}%
}}}}
         \put(9600,-4200){\makebox(0,0)[lb]{\smash{{\SetFigFont{12}{14.4}{\familydefault}{\mddefault}{\updefault}{\color[rgb]{0,0,0}$2132^a$}%
}}}}

         \put(1850,-3200){\makebox(0,0)[lb]{\smash{{\SetFigFont{12}{14.4}{\familydefault}{\mddefault}{\updefault}{\color[rgb]{0,0,0}$13123^e$}%
}}}}
         \put(3050,-3200){\makebox(0,0)[lb]{\smash{{\SetFigFont{12}{14.4}{\familydefault}{\mddefault}{\updefault}{\color[rgb]{0,0,0}$12321^c$}%
}}}}
         \put(4200,-3200){\makebox(0,0)[lb]{\smash{{\SetFigFont{12}{14.4}{\familydefault}{\mddefault}{\updefault}{\color[rgb]{0,0,0}$12132^a$}%
}}}}
         \put(5450,-3200){\makebox(0,0)[lb]{\smash{{\SetFigFont{12}{14.4}{\familydefault}{\mddefault}{\updefault}{\color[rgb]{0,0,0}$31232^e$}%
}}}}
         \put(6550,-3200){\makebox(0,0)[lb]{\smash{{\SetFigFont{12}{14.4}{\familydefault}{\mddefault}{\updefault}{\color[rgb]{0,0,0}$23213^d$}%
}}}}
         \put(7900,-3200){\makebox(0,0)[lb]{\smash{{\SetFigFont{12}{14.4}{\familydefault}{\mddefault}{\updefault}{\color[rgb]{0,0,0}$21323^a$}%
}}}}
         \put(10100,-3200){\makebox(0,0)[lb]{\smash{{\SetFigFont{12}{14.4}{\familydefault}{\mddefault}{\updefault}{\color[rgb]{0,0,0}$21321^b$}%
}}}}

         \put(2300,-1700){\makebox(0,0)[lb]{\smash{{\SetFigFont{12}{14.4}{\familydefault}{\mddefault}{\updefault}{\color[rgb]{0,0,0}$123213^d$}%
}}}}
         \put(4600,-1700){\makebox(0,0)[lb]{\smash{{\SetFigFont{12}{14.4}{\familydefault}{\mddefault}{\updefault}{\color[rgb]{0,0,0}$131232^e$}%
}}}}
         \put(5850,-1900){\makebox(0,0)[lb]{\smash{{\SetFigFont{12}{14.4}{\familydefault}{\mddefault}{\updefault}{\color[rgb]{0,0,0}$121321^b$}%
}}}}
         \put(7100,-1700){\makebox(0,0)[lb]{\smash{{\SetFigFont{12}{14.4}{\familydefault}{\mddefault}{\updefault}{\color[rgb]{0,0,0}$121323^a$}%
}}}}
         \put(8250,-2000){\makebox(0,0)[lb]{\smash{{\SetFigFont{12}{14.4}{\familydefault}{\mddefault}{\updefault}{\color[rgb]{0,0,0}$232132^d$}%
}}}}
         \put(9500,-1900){\makebox(0,0)[lb]{\smash{{\SetFigFont{12}{14.4}{\familydefault}{\mddefault}{\updefault}{\color[rgb]{0,0,0}$213213^b$}%
}}}}

         \put(3400,-600){\makebox(0,0)[lb]{\smash{{\SetFigFont{12}{14.4}{\familydefault}{\mddefault}{\updefault}{\color[rgb]{0,0,0}$1232132^d$}%
}}}}
         \put(5200,-600){\makebox(0,0)[lb]{\smash{{\SetFigFont{12}{14.4}{\familydefault}{\mddefault}{\updefault}{\color[rgb]{0,0,0}$1213213^b$}%
}}}}
         \put(8300,-600){\makebox(0,0)[lb]{\smash{{\SetFigFont{12}{14.4}{\familydefault}{\mddefault}{\updefault}{\color[rgb]{0,0,0}$2132132^c$}%
}}}}

         \put(5800,800){\makebox(0,0)[lb]{\smash{{\SetFigFont{12}{14.4}{\familydefault}{\mddefault}{\updefault}{\color[rgb]{0,0,0}$12132132^c$}%
}}}}
       \end{picture}%
       \caption{The weak order for Example \ref{ex:bruhat} in object c}
       \label{fi:bruhatc}
    \end{figure}

    The rank of the Cartan scheme is three and the length of the longest element
    in $\Homsto a$ (see below) is $8$,
    and hence none of the posets $\Homsto a$, $\Homsto b$ and $\Homsto c$
    with the weak order depicted 
    in Figures \ref{fi:bruhata}, \ref{fi:bruhatb} and \ref{fi:bruhatc}
    can be obtained from a Coxeter group.
    In this respect a particularly interesting case is
    Figure \ref{fi:bruhatc}. Note that for Coxeter groups $W$
    the polynomial $\sum_{w \in W} t^{\ell(w)}$ is a product of factors of the
    form $1+t+\cdots +t^e$. 
    In particular it follows that the coefficient sequence of $\sum_{w \in W} t^{\ell(w)}$ is unimodal, i.e.,
    weakly increases and weakly decreases along increasing $t$ powers.
    Now despite the fact that they cannot arise from Coxeter groups for Figure \ref{fi:bruhata} and \ref{fi:bruhatb} 
    the analogously defined polynomial still has the nice factorization. But in the example Figure \ref{fi:bruhatc} 
    this fails and moreover the coefficient sequence $1,3,6,7,6,7,6,3,1$ is not unimodal.
  \end{example}

  In what follows, for all $a\in A$
  we consider $\Homsto a$ as a poset with respect to the weak
  order.

  \begin{lemma} \label{le:maxchains}
    Let $a\in A$. Then
    all maximal chains in $\Homsto a$ have the same length. This number is
    independent of $a$ in the connected component of $\cC $ containing $a$.
    Hence, $\Delta( \Homsto a )$ is a pure simplicial complex.
  \end{lemma}

  \begin{proof}
    A chain $u_0\swo u_1\swo \cdots \swo u_k$ in $\Homsto a$,
    where $k\in \ndN _0$,
    is maximal if and only if
    $\ell (u_j)=j$ for all $j\in \{0,1,\dots,k\}$ and
    \begin{align}
        \ell (u_k\s _i)\le \ell (u_k)\quad \text{for all $i\in I$.}
        \label{eq:maxch}
    \end{align}
    Lemma~\ref{le:posi} and \eqref{eq:maxch} imply that
    $u_k(\al _i)\in -R^a_+$ for all $i\in I$. Hence $u_k(\al )\in -R^a_+$
    for all $\al \in R^b_+$, where $b\in A$ such that $u_k\in \Hom (b,a)$.
    Then $k=\ell (u_k)=|R^b_+|=|R^a_+|=|R^a|/2$ by Lemma~\ref{le:length}.
    In the connected component of $\cC $ containing $a$
    the number of roots per object is constant by Axiom (R3).
  \end{proof}

  \begin{corollary} \label{co:minimax}
    Let $a\in A$ and $J\subseteq I$.
    There is a unique minimal and a unique maximal element in
    $\Homsto[\Wg _J(\cC )]a$.
  \end{corollary}

  \begin{proof}
    By Proposition~\ref{pr:Wgfun} the groupoid $\Wg _J(\cC )$ is isomorphic to
    the Weyl groupoid of a Cartan scheme. The length function on
    $\Wg _J(\cC )$ is $\ell _J$, which itself coincides with
    the restriction of the length function
    of $\Wg (\cC )$ by Proposition~\ref{co:lJ=l}. Thus we may assume that
    $J=I$.

    The unique minimal element in $\Homsto a$ is $\id^a$.
    In view of the proof of
    Lemma~\ref{le:maxchains}, maximal elements have length $|R^a_+|$.
    By \cite[Cor.\,5]{a-HeckYam08} there is a unique element in $\Homsto a$ of
    maximal length, which implies the claim.
  \end{proof}

  \begin{definition} \label{de:maxel}
    For all $a\in A$ and $J\subseteq I$ we write $w_J$ for the unique maximal
    element of $\Homsto[\Wg _J(\cC )]a$ with respect to weak order. We say
    that $w_J$ is the {\it longest word} over $J$. 
  \end{definition}

  The element $w_J$ in Definition~\ref{de:maxel} depends on the object $a$.
  Nevertheless for brevity we omit $a$ in the notation, since usually it is
  clear from the context what it is.

  \begin{lemma} \label{le:contentmaximal}
    Let $a\in A$, $J\subseteq I$ and $w_J$ the unique maximal
    element of $\Homsto[\Wg _J(\cC )]a$ with respect to weak order. 
    Then $J(w_J) = J$. 
  \end{lemma}

  \begin{proof} 
     This follows from Lemma \ref{le:posi}.
  \end{proof}

  In \cite[p.\,17]{b-BjoeBre05} left descent sets and left descents of
  elements of Coxeter groups have been defined. We generalize the
  definition to our setting, and introduce a related notion.

  For all $a\in A$ and $w\in \Homsto a$ let
  \begin{align}
    D_L(w)=&\{ s\in \Homsto a\,|\,\ell (s)=1,\,s\wo w\},\\
    I_L(w)=&\{i\in I\,|\,\id^a\s _i \in D_L(w)\}.
    \label{eq:DL}
  \end{align}
  The set $D _L(w)$ is called the \textit{left
  descent set} of $w$ and its elements are called the left descents of $w$.
  Clearly, every element $w\not=\id^a$ has left descents.
  Similarly, let
  \begin{align}
    \bar{D}_L(w)=\{ w_J\in \Homsto a\,|\,J\subseteq I,\,w_J\wo w\}.
    \label{eq:barDL}
  \end{align}
  Since $w_{\{j\}}=\id^a \s _j$ for all $j\in I$, we have a natural inclusion
  $D_L(w)\subseteq \bar{D}_L(w)$.
  In the sequel we will consider $\bar{D}_L(w)$ as a subposet of $\Homsto a$
  ordered by the weak order.

 \begin{lemma} \label{le:welldefined}
    Let $a \in A$, $w \in \Homsto a$ and $J = I_L(w)$. Then $w_J \wo w$. 
 \end{lemma}

 \begin{proof}
    Set $x:=w^{-1}w_J$. Then $w=w_J x^{-1}$.
    To prove that $w_J\in \bar{D}_L(w)$, we have to show that $\ell (x)=\ell
    (w)-\ell (w_J)$. By definition of $I_L(w)$ and Lemma~\ref{le:posi}
    we conclude that $w^{-1}(\al _j)\in -\ndN _0^I$ and $w_J(\al _j)\in
    -\lspan _{\ndN _0}\{\al _m\,|\,m\in J\}$ for all $j\in J$.
    Hence $x(\al _j)\in \ndN _0^I$ for all $j\in J$.
    Therefore $x\in \Wg ^J(\cC )$ by Lemma~\ref{le:posi}, and hence
    $\ell (x w_J^{-1})=\ell (x)+\ell (w_J^{-1})$ by
    Proposition~\ref{pr:decomp}(2). This yields the claim.
 \end{proof}

 \begin{proposition} \label{pr:extdesc}
    Let $a\in A$ and $w\in \Homsto a$.
    The map
    $2^{I_L(w)} \to \bar{D}_L(w)$, $J\mapsto w_J$, is an isomorphism
    of posets.
  \end{proposition}

  \begin{proof}
    \noindent {\sf Well-defined:}
    By Lemma \ref{le:welldefined} the map $2^{I_L(w)}\to \bar{D}_L(w)$ is well
    defined.

    \noindent {\sf Injectivity:}
    This follows immediately from Lemma~\ref{le:contentmaximal}.

    \noindent {\sf Surjectivity:}
    Let $J\subseteq I$ such that $w_J\wo w$.
    The definition of $w_J$
    implies that $\id ^a\s _j\wo w_J$ for all $j\in J$, and hence $J\subseteq
    I_L(w_J)\subseteq I_L(w)$.
    Thus the map $2^{I_L(w)} \to \bar{D}_L(w)$
    is surjective.

    \noindent {\sf Poset-Isomorphism:}
    Definition~\ref{de:maxel} implies that
    $w_J\wo w_{J'}$ whenever $J\subseteq J'\subseteq I$. Conversely, let
    $J,J'\subseteq I$ with $w_J\wo w_{J'}$. By Corollary
    \ref{le:contentmaximal} it follows that $J = J(w_{J})$ and 
    $J' = J(w_{J'})$. Hence from $w_J\wo w_{J'}$ we infer $J \subseteq J'$.
  \end{proof}

  \begin{proposition} \label{pr:intervaliso}
    Let $a,b\in A$, $u\in \Hom (b,a)$ and $v\in \Homsto a$ such that $u\swo
    v$.

    \begin{itemize}
       \item[(1)] The map $w\mapsto u^{-1}w$ is an isomorphism of posets
          from the interval $[u,v]$ to the interval $[\id^b,u^{-1}v]$.

       \item[(2)] The map $w\mapsto u^{-1}w$ is an isomorphism of posets
          from the interval $(u,v)$ to the interval $(\id^b,u^{-1}v)$.
    \end{itemize}
  \end{proposition}

  \begin{proof} Follow the proof of \cite[Prop.\,3.1.6]{b-BjoeBre05}.
    This uses only basic properties of the length function which hold also
    for the length function of $\Wg (\cC )$.
    The arguments are the same for both (1) and (2), and work also if one
    considers intervals which are open on one side and closed on the other.
  \end{proof}

  Let $(P,\le )$ be a poset and $U\subseteq P$ a subset. An element $z\in P$
  is called the \textit{meet} of $U$ if
  \begin{itemize}
    \item $z\le u$ for all $u\in U$, and
    \item $y\le z$ for all $y\in P$ with $y\le u$ for all $u\in U$.
  \end{itemize}
  If it exists, the meet of $U$ is unique and is denoted by $\bigwedge U$.
  The meet of two elements $x,y\in P$ is denoted by $x\wedge y$.
  Similarly, an element $z\in P$
  is called the \textit{join} of $U$ if
  \begin{itemize}
    \item $u\le z$ for all $u\in U$, and
    \item $z\le y$ for all $y\in P$ with $u\le y$ for all $u\in U$.
  \end{itemize}
  If it exists, the join of $U$ is unique and is denoted by $\bigvee U$.
  The join of two elements $x,y\in P$ is denoted by $x\vee y$.
  In the sequel we write $\vee$ for the join 
  and $\wedge$ for the meet in $\Homsto a$
  with respect to the weak order.

  A poset is called a \textit{meet semilattice}, if every finite non-empty
  subset has a meet. Finite Coxeter groups with weak order form a meet
  semilattice by \cite[Thm.\,3.2.1]{b-BjoeBre05}, but the proof uses the
  exchange condition which is not available in our setting (see Remark 
  \ref{re:cmsl} for the case of infinite Coxeter groups and Weyl groupoids). 
  We present for Weyl groupoids of Cartan schemes
  a proof which is based on Proposition~\ref{pr:extdesc}.
  The following lemma is one step in our proof.

  \begin{lemma}
    Let $a\in A$ and $u,v,w\in \Homsto a$ such that $w\wo u$ and $w\wo v$.
    If $I_L(w)\subsetneq I_L(u)\cap I_L(v)$ then
    there exists $w'\in \Homsto a$ such that $w\swo w'$ and $w'\wo u$,
    $w'\wo v$.
    \label{le:notmax}
  \end{lemma}

  \begin{proof}
    We proceed by induction on the length of $w$.
    If $\ell (w)=0$ then $w=\id^a$ and the claim holds with $w'=w_{I_L(u)\cap
    I_L(v)}$ by Lemma \ref{le:welldefined}.

    Assume now that $\ell (w)>0$. Let $J=I_L(u)\cap I_L(v)$, and let
    $w_0\in \Homsto[\Wg _J(\cC )]{a}$ be maximal with respect to weak order
    such that $w_0\wo w$. Then $\ell (w_0)>0$ since $\ell (w)>0$ and
    $I_L(w) \subseteq J$.
    Further, $w_0\not=\id^a w_J$ since $I_L(w)\not=J$.
    Let $b\in A$ and $u_1,v_1,w_1\in \Homsto b$ such that $w=w_0w_1$,
    $u=w_0u_1$, and $v=w_0v_1$.
    Then $w_0\wo u$ and $w_0\wo v$ by transitivity of $\wo $, and hence
    $w_1\wo u_1$, $w_1\wo v_1$ by Proposition~\ref{pr:intervaliso}.
    Moreover, $I_L(w_1)\cap J=\emptyset $ by the maximality
    of $w_0$, and $I_L(u_1)\cap I_L(v_1)\cap J\not=\emptyset $ since
    $w_0\not=\id^a w_J$. Since $\ell (w_1)<\ell (w)$, induction hypothesis
    provides us with $w''\in \Homsto b$ such that $w_1\swo w''$ and $w''\wo
    u_1$, $w''\wo v_1$. Then the lemma holds with $w'=w_0w''$ by
    Proposition~\ref{pr:intervaliso}.
  \end{proof}

  \begin{theorem} \label{th:semilattice}
    Let $a\in A$. Then $\Homsto a$ is a meet semilattice.
  \end{theorem}

  \begin{proof}
    For all $v\in \Homsto a$ the set $\{w\in \Homsto a\,|\,w\wo v\}$ is
    finite. Hence it suffices to show that any pair of elements of $\Homsto a$
    has a meet.

    Let $u,v\in \Homsto a$.
    We prove by induction
    on the length of $u$ that the set $\{u,v\}$ has a meet.

    For all $w\in \Homsto a$ with $w\wo u$ and $w\wo v$
    it follows that $I_L(w)\subseteq I_L(u)\cap I_L(v)$.
    Thus if $I_L(u)\cap I_L(v)=\emptyset $, then $w=\id^a$, and hence
    $u\wedge v=\id^a$. This happens in particular if $\ell (u)=0$.

    Assume now that $J:=I_L(u)\cap I_L(v)\not=\emptyset $, and let $w_1,w_2\in
    \Homsto a$ be maximal with respect to weak order
    such that $w_i\wo u$ and $w_i\wo v$ for all $i\in \{1,2\}$. We show that
    $w_1=w_2$. The maximality assumption and Lemma~\ref{le:notmax} imply that
    $I_L(w_1)=I_L(w_2)=J$. Hence $\id^a w_J\wo w_i$ for all $i\in \{1,2\}$
    by Lemma \ref{le:welldefined}. Therefore there exist unique $b\in A$,
    $u',v',w'_1,w'_2\in \Homsto b$ such that $\id^a w_J\in \Hom (b,a)$,
    $w_i=\id^a w_J w'_i$, $u=\id^a w_J u'$, $v=\id^a w_J v'$.
    Proposition~\ref{pr:intervaliso} implies that $w'_1,w'_2$ are maximal.
    Since $\ell (u')<\ell (u)$, induction hypothesis implies that $w'_1=w'_2$,
    and hence $w_1=w_2$. Thus the theorem is proven.
  \end{proof}

  \begin{remark} \label{re:cmsl}
    The proof of Theorem \ref{th:semilattice} does not use the assumption that $\Homsto a$ is
    finite. Thus analogously to the case of Coxeter groups (see
    \cite[Thm.\,3.2.1]{b-BjoeBre05})
    in the weak order of Weyl groupoids the meet of an arbitrary
    subset exists and therefore the weak order forms a complete meet semilattice.
  \end{remark}

  Since $\Homsto a$ is finite and has a unique maximal element by
  Corollary~\ref{co:minimax}, the following corollary follows from
  Theorem~\ref{th:semilattice}
  by standard arguments in lattice theory.

  \begin{corollary} \label{co:lattice}
    Let $a\in A$. Then $\Homsto a$ is a lattice.
  \end{corollary}
  
  The following result is the extension to Weyl groupoids of Theorem 3.2.7 from
  \cite{b-BjoeBre05}.

  \begin{theorem} \label{th:interval}
     Let $a \in A$ and $u,v \in \Homsto{a}$ such that $u \wo v$.
     Let $J=I_L(u^{-1}v)$.
     If $u^{-1} v \not= w_J$ then $(u,v)$ is contractible.
     If $u^{-1} v = w_J$ then $(u,v)$ is
     homotopy equivalent to a sphere of dimension $|J|-2$. 
  \end{theorem}

  \begin{proof}
     By Proposition \ref{pr:intervaliso} it follows that we only need to
     consider the case $u = \id^a$. 
     Consider the map 
     $f : (\id^a,v) \rightarrow (\id^a,v)$ sending $w \in 
     (\id^a,v)$ to $w_{I_L(w)}$. 

     Let $w,w' \in \Homsto a$ with $w \wo w'$. Then $I_L(w) \subseteq I_L(w')$
     and hence $f(w) \wo f(w') \wo w'$. 
     Hence by Theorem \ref{th:closure} it follows that $(\id^a,v)$ and its image under $f$
     are homotopy equivalent.
     From Proposition \ref{pr:extdesc} we infer that the image of 
     $[\id^a,v]$ under $f$ is as a poset isomorphic to $2^{I_L(v)}$
     ordered by inclusion. 

     If $v = w_{I_L(v)}$ then 
     Proposition \ref{pr:extdesc} implies that the image of the
     open interval $(\id^a,v)$ under $f$ is isomorphic to the open 
     interval $(\emptyset, I_L(v))$ and hence by Example \ref{ex:boolean}
     homeomorphic to a $|I_L(v)|-2$ sphere. If $v \neq w_{I_L(v)}$ then 
     $w_{I_L(v)}$ is the unique maximal element of the image of $(\id^a,v)$ 
     under $f$. In particular, the image is isomorphic to the half open
     interval $(\emptyset, I_L(v)]$. Since a poset with unique maximal
     element is contractible the rest of the assertion follows.
  \end{proof}

  \begin{remark}
    \label{re:involution}
    For all $a\in A$ let $\tau (a)\in A$ such that $w_I\in \Hom (\tau (a),a)$.
    Since $w_I$ maps positive roots to negative roots,
    Lemma~\ref{le:posi} implies that $w_I^{-1}$ is a maximal element in
    $\Hom (a,\tau (a))$. Hence $\tau ^2(a)=a$ by Corollary~\ref{co:minimax} and
    the definition of $\tau $. Thus $\tau :A\to A$, $a\mapsto \tau (a)$,
    is an involution of $A$.
  \end{remark}

  The longest element of a Weyl group induces an automorphism of the
  corresponding Dynkin diagram. This automorphism can be generalized to Weyl
  groupoids as follows.
  Let $a\in A$. Since $w_I\in \Hom (a,\tau (a))$ maps positive roots to
  negative roots, Axiom (R1) implies that there exists a permutation
  $\tau _I^a\in \Perm I$ such that $w_I\id^a(\al _j)=-\al _{\tau _I^a(j)}$.

  \begin{lemma}\ \label{le:tauI}
    \begin{itemize}
       \item[(1)] For all $a\in A$ the permutation $\tau _I^a$ is an involution and
          $\tau _I^b=\tau _I^a$ for all $b\in A$ in the connected component of $a$
          in $\cC $.
       \item[(2)] For all $a\in A$ and $i\in I$ we have $w_I\s _i^a w_I=\s _{\tau_I^a(i)}^{\tau (a)}$.
    \end{itemize}
  \end{lemma}

  \begin{proof}
    The definition of $\tau _I^a$ and the formula $w_I w_I\id^a=\id^a$
    imply that $\tau _I^{\tau (a)}\tau _I^a=\id $ for all $a\in A$.
   
    \begin{itemize}
      \item[(2)] Let $a\in A$ and $i,j\in I$. Then
         $w_I\s _i w_I\s _j^{\rfl _j(\tau (a))}
         \in \Hom (\rfl _j(\tau (a)),\tau (\rfl _i(a)))$.
         Assume that $\tau _I^{\tau (a)}(j)=i$, that is, $j=\tau _I^a(i)$.
         Then
         \begin{align} \label{eq:wsws}
            w_I\s _i w_I\s _j^{\rfl _j(\tau (a))}(\al _j)
           =-w_I\s _i w_I\id^{\tau (a)}(\al _j)
           =w_I\s _i^a(\al _i)
           =-w_I\id^{\rfl _i(a)}(\al _i)
           =\al _{\tau _I^{\rfl _i(a)}(i)}.
         \end{align}
         Moreover, $w_I\s _i w_I\s _j^{\rfl _j(\tau (a))}$ maps any positive root
         different from $\al _j$ to a positive root since $w_I$ maps positive roots
         to negative roots and for all $l\in I$, $b\in A$ the map $\s _l^b$
         sends positive roots different from $\al _l$ to positive roots, see
         \cite[Lemma\,1]{a-HeckYam08}.
         Thus $\ell (w_I\s _i w_I\s _j^{\rfl _j(\tau (a))})=0$ by
         Lemma~\ref{le:length}, and hence
         $w_I\s _i^a w_I=\s _j^{\tau (a)}$.
    
      \item[(1)] Since for all $a\in A$ the object $\tau (a)$
         is in the same connected component as $a$,
         it suffices to show that for all $a\in A$ and $i\in I$ the permutations
         $\tau _I^a$ and $\tau _I^{\rfl _i(a)}$ are equal.
         Let $a\in A$ and $i\in I$. By (2) we obtain that
         \begin{align}
            \s _{\tau _I^a(i)}w_I\s _i^a=w_I\id^a,
            \label{eq:wsw=s}
         \end{align}
         and Equation~\eqref{eq:wsws} gives that $\tau _I^{\rfl _i(a)}(i)=j=\tau_I^a(i)$. 
         Applying Equation~\eqref{eq:wsw=s} to all $\al _k$
         with $k\in I$ implies that $\tau _I^a=\tau _I^{\rfl _i(a)}$.
    \end{itemize}
  \end{proof}

  For all $a\in A$ define the map
  $t^a:\Homsto a\to \Homsto{\tau (a)}$ by
  \[ t^a(\id^a \s _{i_1}\cdots \s _{i_k})=
  \id^{\tau (a)}\s _{\tau _I^a(i_1)}\cdots \s _{\tau _I^a(i_k)}\quad
  \text{for all $k\in \ndN _0$, $i_1,\dots,i_k\in I$.} \]

%
%

  \begin{proposition}
    \label{pr:posetisomorphism}
    Let $a\in A$. Then $t^a(w)=w_I w w_I$ and $\ell (t^a(w))=\ell (w)$
    for all $w\in \Homsto a$.
    The map $t^a$
    is an isomorphism of posets with respect to weak order.
  \end{proposition}

  \begin{proof}
    Lemma~\ref{le:tauI}(1) and (2) imply that
    \begin{align*}
      \id^{\tau (a)}w_I\s _{i_1}\cdots \s _{i_k}w_I
      =&\id^{\tau (a)}(w_I\s _{i_1}w_I)(w_I\s _{i_2}w_I)\cdots (w_I\s
      _{i_k}w_I)\\
      =&\id^{\tau (a)}\s _{\tau _I^a(i_1)}\s _{\tau _I^a(i_2)}\cdots
      \s _{\tau _I^a(i_k)}
    \end{align*}
    for all $a\in A$, $k\in \ndN _0$, and $i_1,\dots,i_k\in I$.
    Hence
    $t^a$ is well-defined and
    the first claim holds.
    Since $w_I w_I\id^a=\id^a$, we conclude that
    $t^{\tau (a)}t^a (w)=w$ for all $w\in \Homsto a$
    and
    $t^a t^{\tau (a)}(w)=w$ for all $w\in \Homsto{\tau (a)}$, and hence $t^a$
    is bijective. It is clear from the definition and bijectivity of $t^a$
    that $t^a$ preserves length and therefore it preserves
    and reflects weak order.
  \end{proof}

  A lattice $P$ 
  with unique minimal element $\hat{0}$ and unique maximal element $\hat{1}$ 
  is called \textit{ortho-complemented} if there is a map
  $\perp : P \rightarrow P$ such that (O1) $p \wedge p^\perp = \hat{0}$, (O2)
  $p \vee p^\perp = \hat{1}$, (O3) For all $p \in P$ we have $(p^\perp)^\perp
  = p$ and (O4) for all $p \preceq q$ in $P$ we have $q^\perp \preceq
  p^\perp$.

  \begin{lemma} \label{le:complement}
     Let $a \in A$ and $w \in \Homsto a$. Then the following hold.
     \begin{itemize}
        \item[(1)] $\ell(w) + \ell(ww_I) = \ell(w_I)$. 
        \item[(2)] $I_L(w) \cap I_L(ww_I) = \emptyset$.
        \item[(3)] For $i \in I$ we have $i \in I_L(w)$ if 
                 and only if $i \not\in I_L(ww_I)$.
     \end{itemize}
  \end{lemma} 
  \begin{proof}\ 
    \begin{itemize}
       \item[(1)] For any $b\in A$ and $v \in \Hom (b,a)$ we have 
          $\ell(v) = \# \{ \alpha \in R^b_+\,|\,v(\alpha) \in -R^a_+ \}$.
          Now $w_I (\alpha) \in -R^b_+$ for all $\alpha \in R^{\tau (b)}_+$.
          Thus for $\alpha \in R^b_+$ we have 
          $$w(\alpha) \in -R^a_+ \Leftrightarrow
          w w_I(-w_I(\alpha)) \in R^a_+.$$
          This implies that $\ell(w) + \ell(ww_I) = \ell(w_I)$.
       \item[(2)] Let $i \in I_L(w) \cap I_L(ww_I)$.  
          Then $\ell(\s_iw) = \ell(w) -1$ and
          $\ell(\s_i ww_I) = \ell(ww_I) -1$.
          Hence 
          \begin{eqnarray*} 
            \ell(\s_i w) + \ell(\s_i ww_I) & = & \ell(w) -1 + \ell(ww_I) -1 \\
                                           & = &  \ell(w_I)-2.
          \end{eqnarray*}
          This contradicts (1) and hence $I_L(w) \cap I_L(ww_I) = \emptyset$.
       \item[(3)] By (2) it suffices to show that $I_L(w) \cup I_L(ww_I) = I$.
          Assume there is an $i \in I \setminus (I_L(w) \cup I_L(ww_I))$.
          Then $\ell(\s_iw) = \ell(w) +1$ and $\ell(\s_i ww_I) = \ell(ww_I)+1$.
          Analogously to (2) we obtain that
          \begin{eqnarray*} 
            \ell(\s_i w) + \ell(\s_i ww_I) & = & \ell(w) +1 + \ell(ww_I) +1 \\
                                           & = &  \ell(w_I)+2
          \end{eqnarray*}
          which is a contradiction to (1) and we are done.
    \end{itemize}
  \end{proof}

  \begin{theorem} \label{th:orthocomplemented}
     Let $a \in A$. Then the map $\perp : \Homsto a \rightarrow \Homsto a$
     defined by $w^\perp := w w_I$ satisfies (O1) - (O4). Thus 
     $\Homsto{a}$ with the weak order is an ortho-complemented 
     lattice.
  \end{theorem}

  \begin{proof}
     \begin{itemize}
        \item[(O1)] This follows immediately from Lemma \ref{le:complement} (2).
        \item[(O2)]
           By Lemma~\ref{le:complement}(3) we know that
           $I_L(w) \cup I_L(ww_I) = I$. Thus any $v\in \Homsto a$
           with $w\wo v$, $ww_I\wo v$ satisfies $w_I\wo v$ by
           Lemma~\ref{le:welldefined}. Hence $w\vee ww_I=w_I$.
        \item[(O3)] This follows from the definition of $\perp$
          and Remark \ref{re:involution}.
        \item[(O4)] Let $u,v\in \Homsto a$ with $u \wo v$.
          If $\ell(u) = 0$ then clearly $v^\perp  \wo u^\perp = w_I$.
          Now proceed by induction on $\ell(u)$. 
          Assume that $\ell(u) \geq 1$ and let $i \in I_L(u)$.
          Then $i \in I_L(v)$ and we find
          $\bar{u}$ and $\bar{v}$ in $\Homsto {\rfl _i(a)}$ 
          such that $u = \s_i \bar{u}$ and $v = \s_i \bar{v}$.
          Then $\bar{u} \wo \bar{v}$. By the induction hypothesis we
          obtain that $\bar{v}^\perp \wo \bar{u}^\perp$. 
          Since $i \not\in I_L(\bar{v})$ and $i \not\in I_L(\bar{u})$
          it follows from Lemma \ref{le:complement} (3) and the definition of 
          $\perp$ that $i \in I_L(\bar{v}^\perp )$ and $i \in
          I_L(\bar{u}^\perp )$. 
          Hence $\s_i \bar{v}^\perp \wo \s_i \bar{u}^\perp$. 
          By the definition of $\perp$
          this implies that $v w_I = \s_i \bar{v}w_I \wo \s_i
            \bar{u}w_I = uw_I$. Hence $v^\perp \wo u^\perp$.
    \end{itemize}
  \end{proof} 
           
  The following proposition strengthens Proposition \ref{pr:extdesc}
  showing that the embedding is indeed an embedding of lattices.
   
  \begin{proposition}
     \label{pr:meetjoin}
     Let $a \in A$ and $J,J' \subseteq I$. 
     Then $w_J \wedge w_{J'} = w_{J \cap J'}$ and $w_J \vee w_{J'} 
     = w_{J \cup J'}$.
     In particular, the map $2^I \rightarrow \Homsto a$, $J \mapsto w_J$ is 
     an embedding of lattices. 
  \end{proposition}
  \begin{proof}
     \begin{itemize}
        \item[($\wedge$)]
         By Proposition \ref{pr:extdesc} it follows that $w_{J \cap J'} \wo
         w_J,w_{J'}$. By Theorem  \ref{th:semilattice} there is a meet
         $w := w_J \wedge w_{J'}$ and hence $w_{J \cap J'} \wo w$.
         Let $b\in A$ such that $w\in \Hom (b,a)$.
         From $w \wo w_J$ and $w \wo w_{J'}$ we deduce that there are
         $u, u' \in \Homsto b$ such that $w_J = w u$, $w_{J'} = w u'$ and 
         $\ell(w_J) = \ell(w) + \ell(u)$, $\ell(w_{J'}) = \ell(w) + \ell(u')$. 
         By $w_{J \cap J'} \wo w$ we deduce that there is
         $v \in \Wg (\cC )$ such that 
         $w = w_{J \cap J'}v$ and $\ell(w) = \ell(w_{J \cap J'}) + \ell(v)$.  
         By $w_{J \cap J'}vu = w_J$ and $w_{J \cap J'}vu' = w_{J'}$ it follows that
         $I_L(v) \subseteq J \cap J'$.
         However, by the fact that $w_{J \cap J'}$ is the longest
         word in $J \cap J'$ and $\ell(w_{J \cap J'}) + \ell(v) = \ell(w_{J \cap J'}v)$
         it follows that $v = \id^a$ and hence $w = w_{J\cap J'}$.
     \item[($\vee$)]
         By Proposition \ref{pr:extdesc} it follows that $w_J,w_{J'} \wo
         w_{J \cup J'}$. Let now $w\in \Homsto a$ such that $w_J,w_{J'}\wo w$.
         We have to show that $w_{J\cup J'}\wo w$.
         By Proposition~\ref{pr:extdesc} with $w=w_J$ we conclude that
         $I_L(w_J)=J$, and similarly $I_L(w_{J'})=J'$.
         Thus $J\cup J'=I_L(w_J)\cup I_L(w_{J'})\subseteq I_L(w)$.
         Lemma~\ref{le:welldefined} and Proposition~\ref{pr:extdesc}
         imply that $w_{J\cup J'}\wo w_{I_L(w)}\wo w$ and we are done.
     \end{itemize}
  \end{proof}

  The following is an immediate consequence of Proposition \ref{pr:meetjoin}.

  \begin{corollary} \label{co:joinD}
    Let $a \in A$. Then for all $J \subseteq I$ we have
    $$\bigvee_{i \in J} \id^a \s_i = \id^a w_{J}.$$
    In particular, for all $w \in \Wg(\cC)$ we have
    $$\bigvee_{i \in I_L(w)} \id^a \s_i = \id^a w_{I_L(w)}.$$
  \end{corollary} 

  Next we present a formula about the factors appearing in a reduced
  decomposition of the meet of two morphisms.
  
  \begin{theorem} \label{th:Ju}
    Let $a\in A$ and $u,v\in \Homsto a$. Then
    \[ J(u)\cup J(v)=J(u\wedge v)\cup J(u^{-1}v). \]
  \end{theorem}
  
  \begin{proof}
    Since $u\wedge v\wo u$ and $u\wedge v\wo v$, it follows that
    $J(u\wedge v)\subseteq J(u)\cap J(v)$. Moreover,
    $J(u^{-1}v)\subseteq J(u)\cup J(v)$,
    and hence the inclusion $\supseteq$ in the theorem holds.

    Now we prove the inclusion $\subseteq$ by induction on $\ell (u)+\ell (v)$.  
    If $\ell (u)=\ell (v)=0$ then the claim clearly holds. Assume now that
    $\ell (u)+\ell (v)>0$.
    
    \textit{Case 1. $u\wedge v\not=\id ^a$.}
    Then there exists $i\in I_L(u)\cap I_L(v)$.
    Let $u_0,v_0\in \Homsto {\rfl _i(a)}$ such that
    $u=\s _i u_0$, $v=\s _i v_0$.
    Then
    \begin{align*}
      J(u)=J(u_0)\cup \{i\},\, J(v)=J(v_0)\cup \{i\},\,
      J(u\wedge v)=J(\s _i(u_0\wedge v_0))=J(u_0\wedge v_0)\cup \{i\},
    \end{align*}
    and
    $u^{-1}v=u_0^{-1}v_0$. Thus the claim follows from the induction hypothesis.
    
    \textit{Case 2. $u\wedge v=\id ^a$, $J(u)\not\subseteq J(v)$.}
    By Proposition~\ref{pr:decomp}
    there exist unique elements
    $u^J\in \Wg ^{J(v)}(\cC )$, $u_J\in \Wg _{J(v)}(\cC )$
    such that $u^{-1}=u^J u_J$.
    Then $u=u_J^{-1}(u^J)^{-1}$ and $\ell (u^J)+\ell (u_J)=\ell (u)$
    and hence $J(u^J)\cup J(u_J)=J(u)$.
    We have $\ell (u_J)<\ell (u)$ since $u\notin \Wg _{J(v)}(\cC )$.
    Further,
    \begin{align} \label{eq:wJ-1}
      u_J^{-1}\wedge v=\id ^a
    \end{align}
    since $u\wedge v=\id ^a$.
    Thus
    \begin{align*}
      J(u)\cup J(v)=&J(u^J)\cup J(u_J^{-1})\cup J(v)
      =J(u^J)\cup (J(u_J^{-1}\wedge v)\cup J(u_J v))\\
      =& J(u^J)\cup J(u_J v)=J(u^J u_J v)=J(u^{-1} v).
    \end{align*}
    Here the second equation holds by induction hypothesis and
    the third by Equation~\eqref{eq:wJ-1}. The fourth equation follows from
    $u_J v\in \Wg _{J(v)}(\cC )$, $u^J\in \Wg ^{J(v)}(\cC )$ and
    Proposition~\ref{pr:decomp}(2).
    
    \textit{Case 3. $u\wedge v=\id ^a$, $J(v)\not\subseteq J(u)$.}
    Replace $u$ and $v$ and apply Case 2.
    
    \textit{Case 4. $u\wedge v=\id ^a$, $J(u)=J(v)$.} Let $J=J(u^{-1}v)$.
    We have to show that $J(u)\subseteq J$.
    By Corollary~\ref{co:minlenrep} there exists a unique minimal element
    $w\in u\Wg _J(\cC )$.
    Since $v=u (u^{-1}v)$ and $u^{-1}v\in \Wg _J(\cC )$,
    there exist $u_1,v_1\in \Wg _J(\cC )$ such that
    \[ u=wu_1,\quad v=wv_1,\quad \ell (u)=\ell (w)+\ell (u_1),
    \quad \ell (v)=\ell (w)+\ell (v_1). \]
    Therefore $w\le u\wedge v=\id ^a$, and hence
    $u\in w\Wg _J(\cC )=\Wg _J(\cC )$.
    Thus $J(u)\subseteq J$.
  \end{proof}

  \section{Coxeter Complex}
  \label{sec:coxeter}

  Throughout this section let
  $\cC =\cC (I,A,(\rfl _i)_{i\in I},(\Cm ^a)_{a\in A})$ be a
  Cartan scheme and let $a \in A$.
  Assume that $\rsC \re (\cC )$ is a finite root system of type $\cC $.

  \begin{definition} \label{de:coxcomp}
    Let $$\Omega_\cC^a 
    := \{ w\Wg_J(\cC)~|~w \in \Homsto{a}, J \subseteq I, |J| = |I|-1 \}.$$
    We call the subset $\Delta_\cC^a$ of the powerset $2^{\Omega_\cC^a}$ whose
    elements are the non-empty subsets
    $F \subseteq \Omega_\cC^a$ such that
    $$\bigcap_{w\Wg_J(\cC) \in F } w\Wg_J(\cC) \neq \emptyset$$
    the {\em Coxeter complex} of $\cC$ at $a$.
  \end{definition}

  By definition, the Coxeter complex  $\Delta_\cC^a$ is a simplicial complex.
  Note that for technical reasons our simplicial complexes do not contain the
  empty set.
  Our goal in this section will be to give a second construction of the
  Coxeter complex. This way we obtain additional information on the structure
  of faces.

  \begin{lemma} \label{le:cosetint}
     Let $J,K\subseteq I$
     and $u,v\in \Homsto a$ such that $u\Wg _J(\cC ) \cap v\Wg _K(\cC )
     \not=\emptyset $. Then
     \begin{align}
       u\Wg _J(\cC )\cap v\Wg _K(\cC )=w\Wg _{J\cap K}(\cC )
     \end{align}
     for some $w\in \Homsto a$.
     In particular, if $J\subseteq K$
     then $w\Wg _J(\cC ) =u\Wg _J(\cC )$ and if
     $J=K$ then $w\Wg _J(\cC ) =u\Wg _J(\cC )=v\Wg _J(\cC )$.
  \end{lemma}

  \begin{proof}
    Assume first that $v=\id ^a$.
    By Proposition \ref{pr:decomp}
    there exist $u_0\in \Wg ^J(\cC )$ and
    $u_1\in \Wg _J(\cC )$ such that $u=u_0 u_1$. Then
    \[ u\Wg _J(\cC )=u_0\Wg _J(\cC ) \]
    and $J(u_0)\subseteq J(w)$
    for all $w\in u\Wg _J(\cC )$ by 
    Corollary \ref{co:minlenrep}.
    Hence $J(u_0)\subseteq J(w)\subseteq K$ for all $w\in
    u\Wg _J(\cC ) \cap v\Wg _K(\cC )$ which is non-empty by assumption.
    Thus
    \begin{align*}
      u\Wg _J(\cC )\cap v\Wg _K(\cC )=
      u_0(\Wg _J(\cC )\cap u_0^{-1}\Wg _K(\cC ))
      =u_0(\Wg _J(\cC )\cap \Wg _K(\cC ))
      =u_0\Wg _{J\cap K}(\cC ).
    \end{align*}
    Let now $v\in \Homsto a$ be an arbitrary element.
    Then
    \begin{align*}
      u\Wg _J(\cC )\cap v\Wg _K(\cC )=
      v(v^{-1}u\Wg _J(\cC )\cap \Wg _K(\cC ))
      =vw_0\Wg _{J\cap K}(\cC )
    \end{align*}
    for some $w_0\in \Homsto b$, where $b\in A$ with $v\in \Hom (b,a)$,
    by the first part of the proof. This implies the claim.
  \end{proof}

  In \cite[Sect.\,1.15]{b-Humphreys90} the Coxeter complex of a reflection
  group was defined by means of hyperplanes in a Euclidean space. We
  introduce an analogous complex for the pair $(\cC ,a)$.
  We show that
  the complex defined this way is isomorphic to the Coxeter
  complex $\Delta _\cC ^a$.

  Let $(\cdot ,\cdot )$ be a scalar product on $\ndR ^I$.
  For any subset $J\subseteq I$ and any $w\in \Homsto a$ let
  \begin{align*}
    \fac^w_J=\{\lambda \in \ndR ^I\,|\,(\lambda ,w(\al _j))=0
    \text{ for all $j\in J$},
    (\lambda ,w(\al _i))>0 \text{ for all $i\in I\setminus J$}\}.
  \end{align*}
  The subsets $\fac ^w_J$ are intersections of hyperplanes and of open
  half-spaces, and are called \textit{faces}.
  For brevity we will omit their
  dependence on the scalar product. 
  By construction the faces do not depend on connected components of $\cC $
  not containing $a$.
  Also, up to the choice of a scalar product the set of faces $\fac^w_J$ 
  does not change when passing from an object $a$ to an object $a'$ from
  a covering Cartan scheme once $a'$ lies in the connected component 
  covering the connected component of $a$.

  The next lemma is the analog of \cite[Lemma\,1.12]{b-Humphreys90}.

  \begin{lemma}
    \label{le:facesformcomplex}
    Let $(\cdot ,\cdot )$ be
    a scalar product on $\ndR ^I$.
    \begin{enumerate}
      \item 
        For all $\lambda \in \ndR ^I$
        there exist $w\in \Homsto a$ and $J\subseteq I$
        such that $\lambda \in \fac ^w_J$.
      \item Let $w_1,w_2\in \Homsto a$ and let $J_1,J_2\subseteq I$.
        If $w_1\Wg _{J_1}(\cC )=w_2\Wg _{J_2}(\cC )$
        then $\fac ^{w_1}_{J_1}=\fac ^{w_2}_{J_2}$.
        If $w_1\Wg _{J_1}(\cC )\not= w_2\Wg _{J_2}(\cC )$
        then $\fac ^{w_1}_{J_1}\cap \fac ^{w_2}_{J_2}=\emptyset $.
    \end{enumerate}
    \label{le:chambers}
  \end{lemma}

  \begin{proof}
    \begin{itemize}
      \item[(1)] Let $k=|\{\beta \in R^a_+\,|\,(\lambda ,\beta )<0\}|$. We
        proceed by
        induction on $k$. If $k=0$ then the claim holds with $w=\id ^a$.

        Assume that $k>0$. Then there exists $i\in I$ such that $(\lambda ,\al
        _i)<0$. Let $\lambda '=\s _i^a (\lambda )$ and define a scalar product
        $(\cdot ,\cdot )'$ on $\ndR ^I$ by
        $(\mu ,\nu )'=(\s _i^{\rfl _i(a)}(\mu ),
        \s _i^{\rfl _i(a)}(\nu ))$ for all $\mu ,\nu \in \ndR ^I$.
        Then
        for all $\beta \in R^{\rfl _i(a)}_+$ we have
        $(\lambda ',\beta )'<0$ if and only if
        $(\lambda ,\s _i^{\rfl _i(a)}(\beta ))<0$. Moreover,
        \[ (\lambda ',\al _i)'=(\s _i^{\rfl _i(a)}\s _i^a(\lambda ),\s 
        _i^{\rfl _i(a)}(\al _i))=-(\lambda ,\al _i)>0, \]
        and $\s _i^{\rfl _i(a)}$ is a bijection between $R^{\rfl
        _i(a)}_+\setminus \{\al _i\}$ and $R^a_+\setminus \{\al _i\}$
        by (R1)--(R3). Hence
        \begin{align*}
          |\{\beta \in R^{\rfl _i(a)}_+\,|\,(\lambda ',\beta )'<0\}|=k-1.
        \end{align*}
        By induction hypothesis there exist $J\subseteq I$
        and $w'\in \Homsto {\rfl _i(a)}$ such that $\lambda '\in \fac ^{w'}_J$.
        Then $\lambda \in \s _i^{\rfl _i(a)}\fac ^{w'}_J=\fac ^w_J$, where
        $w=\s _i^{\rfl _i(a)}w'$.

      \item[(2)] Suppose that $w_1\Wg _{J_1}(\cC )=w_2\Wg _{J_2}(\cC )$.
        Then $J_1=J_2$ and
        $w_2=w_1x$ for some $x\in \Wg _{J_1}(\cC )$. Therefore
        \begin{align*}
          (\lambda ,w_2(\al _i))=(\lambda ,w_1 x(\al _i))=
          (\lambda ,w_1 (\al _i+\sum _{j\in J_1}a_{ji}\al _j))
          =(\lambda ,w_1(\al _i))
        \end{align*}
        for all $\lambda \in \fac ^{w_1}_{J_1}$ and all $i\in I$,
        where $x(\al _i)=\al _i+
        \sum _{j\in J_1}a_{ji}\al _j$ for some $a_{ji}\in \ndZ $ for all $j\in
        J_1$.
        We conclude that $\fac ^{w_1}_{J_1}\subseteq \fac ^{w_2}_{J_2}$, and
        similarly
        $\fac ^{w_2}_{J_2}\subseteq \fac ^{w_1}_{J_1}$. This proves the first
        claim.

        The converse will be proven indirectly.
        Assume that $w_1\Wg _{J_1}(\cC )\not=w_2\Wg _{J_2}(\cC )$ and that
        there exists $\lambda \in \fac ^{w_1}_{J_1}\cap \fac ^{w_2}_{J_2}$.
        Let $b_1,b_2\in A$ such that $w_1\in \Hom (b_1,a)$ and $w_2\in \Hom
        (b_2,a)$. Let $x=w_1^{-1}w_2\in \Hom (b_2,b_1)$.
        By the choice of $\lambda $ and the definition of $x$
        we have $(\lambda ,w_1(\al _j))\ge 0$ and $(\lambda ,w_1 x(\al _j))\ge 0$
        for all $j\in I$. Moreover, equality holds if and only if $j\in J_1$
        respectively $j\in J_2$. 
        Since $x(\al _j)\in R^{b_1}_+\cup -R^{b_1}_+$ for all $j\in I$,
        we conclude that $x(\al _j)\in \sum _{k\in J_1}\ndZ \al _k$ for all $j\in
        J_2$ and that $x(\al _j)\in R^{b_1}_+\setminus \sum _{k\in J_1}\ndZ \al_k$ 
        for all $j\in I\setminus J_2$. Hence $J(x)\subseteq J_1$
        by Lemma~\ref{le:J(w)char}. It follows that
        \begin{align}
          w_2\in w_1\Wg _{J_1}(\cC ).
          \label{eq:w2inw1W}
        \end{align}
        By the first part of the proof we obtain that
        $\fac ^{w_2}_{J_2}=\fac ^{w'_2}_{J_2}$ for all
        $w'_2\in w_2\Wg _{J_2}(\cC )$.
        Hence $J(xx')\subseteq J_1$ for all $x'\in \Wg _{J_2}(\cC )$,
        and therefore $J_2\subseteq J_1$. Symmetry yields that $J_1=J_2$.
        Thus $w_1\Wg _{J_1}(\cC )=w_2\Wg _{J_2}(\cC )$ by
        \eqref{eq:w2inw1W}, a contradiction.
        Hence $F^{w_1}_{J_1}\cap F^{w_2}_{J_2}=\emptyset $.
    \end{itemize}
  \end{proof}

  By definition for any $w \in \Homsto a$ and $J\subseteq I$ 
  the face $\fac ^w_J$ is a relative open polyhedral cone in $\ndR ^I$.
  In particular, it is a relative open cell.
  By Lemma \ref{le:facesformcomplex} the set of all
  $\fac ^w_J$ stratifies $\ndR ^I$. Clearly, this stratification depends on
  the choice of $\cC $, $a$, and the scalar product on $\ndR ^I$.
  In order to show that the stratification indeed gives a regular
  CW-composition of
  $\ndR ^I$ we have to clarify the structure of the closures of the cells.

  \begin{theorem} 
    \label{th:facincl}
    Let $K\subseteq I$ and let $w\in \Homsto a$. 
    Then $\overline{\fac ^w_K}$ is the disjoint union of
    the faces $\fac ^w_J$ for $J \supseteq K$.
    Moreover, for $v\in \Homsto a$ and $J\subseteq I$
    we have $\fac ^v_J\subseteq \overline{\fac ^w_K}$
    if and only if $v\Wg _J(\cC )\supseteq
    w\Wg _K(\cC )$.
  \end{theorem}

  \begin{proof}
    The first assertion follows from the definition of $\fac ^w_K$.

    By Lemma~\ref{le:chambers} the space
    $\ndR ^I$ is the disjoint union of faces, and
    hence $\fac ^v_J\subseteq \overline{\fac ^w_K}$ if and only if
    $\fac ^v_J=\fac ^w_L$ for some $L\supseteq K$.
    Lemma~\ref{le:chambers}(2) implies that the latter is equivalent to
    $v\Wg _J(\cC )=w\Wg _L(\cC )$.
    Clearly, if
    $v\Wg _J(\cC )=w\Wg _L(\cC )$ for some $L\supseteq K$ then
    $v\Wg _J(\cC )\supseteq w\Wg _K(\cC )$.
    Conversely, if
    $v\Wg _J(\cC )\supseteq w\Wg _K(\cC )$ then
    $v^{-1}w\Wg _K(\cC ) \subseteq \Wg _J(\cC )$, and hence
    $v^{-1}w\in \Wg _J(\cC )$ and $K\subseteq J$.
    Thus $v\Wg _J(\cC )=w\Wg _J(\cC )$ and the theorem is proven.
  \end{proof}

  \begin{corollary}
    \label{co:regularcwcomplex}
    The cells 
    $\fac ^w_K$ for $w \in \Homsto a$ and $K \subseteq I$ define a regular
    CW-decomposition of $\ndR ^I$. 
  \end{corollary} 
  \begin{proof}
     From the fact that any root system contains a basis it follows that
     $\fac ^w_I = \{ 0\}$. Hence it follows from Theorem \ref{th:facincl} 
     and the fact that all $\fac^w_J$ are relative open polyhedral cones 
     in $\ndR ^I$ that $\dim \fac^w_J = \# I - \# J$. 
     Since by Lemma \ref{le:facesformcomplex} the cells $\fac^w_J$ are a
     stratification of
     $\ndR ^I$, they actually define a regular CW-decomposition of $\ndR ^I$.
  \end{proof}
 
  Now we define the regular CW-complex $\mathcal{K}^a_{\cC}$ as the 
  regular CW-complex whose cells are the intersections $\fac^w_J \cap 
  S^{\#I - 1}$ of the relative open cones $\fac^w_J$ with the unit
  sphere in $\ndR ^I$ for $J \subseteq I$, $J \neq I$. 
  From Corollary~\ref{co:regularcwcomplex}
  and the fact that all $\fac^w_J$
  are relative open cones with apex in the origin if follows that 
  $\mathcal{K}^a_{\cC}$ is a regular CW-decomposition of $S^{\#I - 1}$.

  \begin{corollary}
    \label{co:regularcwsphere}
    The Coxeter complex $\Delta ^a_{\cC }$ at $a \in A$
    is isomorphic to the complex $\mathcal{K}^a_{\cC}$.
    \label{co:isomorphism}
  \end{corollary}

  \begin{proof}
      Since by Corollary \ref{co:regularcwsphere} the complex
      $\mathcal{K}^a_{\cC}$
      is a regular CW-complex and $\Delta ^a_{\cC }$ is by definition a
      regular CW-complex, it suffices to show that there is an inclusion
      preserving bijection between the faces of $\mathcal{K}^a_{\cC}$ and
      $\Delta ^a_{\cC }$.
 
      Definition~\ref{de:coxcomp} and Lemma~\ref{le:cosetint}
      imply that the faces of the Coxeter complex are in bijection
      with the left cosets $w\Wg _J(\cC )$, where $w\in \Homsto a$
      and $J\subsetneq I$.
      By Lemma~\ref{le:chambers}(2) the faces of $\mathcal{K}^a_{\cC}$
      are also in bijection with these left cosets. Hence it remains to show
      that in both complexes the inclusion of closures of faces corresponds to
      the inclusion of left cosets. For the Coxeter complex
      this holds by definition. For the complex $\mathcal{K}^a_{\cC}$ the claim follows
      from Theorem~\ref{th:facincl}.
  \end{proof}

  Let $\arrg ^a_{\cC}$ be the set of hyperplanes $H_{\alpha} = \{
  \lambda \in \RR^I ~|~(\lambda,\alpha) = 0\}$ for $\alpha \in \rer a_+$.
  Then the complement $\RR^I \setminus \bigcup_{H \in \arrg ^a_{\cC}} H$
  of the 
  arrangement of hyperplanes $\arrg ^a_\cC $ is the
  disjoint union of connected components which are in bijection
  with the maximal 
  faces of $\mathcal{K}_{\cC}^a$.
  It 
  follows by Corollary \ref{co:isomorphism} that $\mathcal{K}_{\cC}^a$ and
  $\Delta ^{a}_{\cC}$ are isomorphic.
  Since $\Delta ^{a}_{\cC }$ is a simplicial complex it follows that all
  connected components of $\RR^I \setminus \bigcup_{H \in \arrg ^a_{\cC }} H$
  are open simplicial cones. 
  In general, an arrangement satisfying this property is called {\it
  simplicial arrangement}.
  
  \begin{corollary}
    The arrangement of hyperplanes $\arrg ^a_\cC $
    is a simplicial arrangement.
  \end{corollary}
  
  From the fact that by Corollary \ref{co:regularcwsphere} the Coxeter complex $\Delta_\cC^a$ is a
  triangulation of a sphere the next corollary follows immediately.

  \begin{corollary}
    \label{co:pure}
    The simplicial complex $\Delta_\cC^a$ is pure of dimension $|I|-1$ and
    each codimension $1$ face of $\Delta_\cC^a$ is contained in exactly two
    faces of maximal dimension. In particular, $\Delta _\cC ^a $ is a
    pseudomanifold.
  \end{corollary}

  Using Theorem \ref{th:facincl}
  one can identify the
  maximal simplices of $\Delta_\cC^a$ with the elements of
  $\Homsto{a}$.
  Hence
  any linear extension of the weak order on $\Homsto a$ defines a
  linear order on the maximal simplices of $\Delta_\cC^a$.
  Indeed it can be shown by the same proof as for the analogous
  statement for Coxeter groups \cite[Theorem 2.1]{a-Bj84} that 
  any linear extension of the weak order defines a shelling order 
  for $\Delta_\cC^a$. The crucial facts about Coxeter groups used by
  Bj\"orner are verified for Weyl groupoids in Lemma \ref{le:cosetint} and 
  Theorem \ref{th:facincl}.

  \begin{theorem} \label{th:shelling}
    Let $\preceq$ be any linear extension of the weak order
    $\leq_R$ on $\Homsto{a}$. Then $\preceq$ is a shelling order
    for $\Delta_\cC^a$.
  \end{theorem}

  We omit the detailed verification of Theorem~\ref{th:shelling}
  here, since the main topological
  consequence Corollary \ref{co:regularcwsphere} is already known. 
  Indeed, 
  Corollary~\ref{co:pure} together with Theorems~\ref{th:shelling} and
  \ref{th:plsphere} imply that $\Delta ^a_\cC $ is a triangulation of a
  PL-sphere.

  \bibliographystyle{amsalpha}
  \bibliography{quantum}
\end{document}